\documentclass{article}
\usepackage{amsmath,amsfonts,enumitem,geometry,mathrsfs}
\usepackage[colorlinks,linkcolor=blue,citecolor=red,urlcolor=blue]{hyperref}
\usepackage[amsmath,thmmarks]{ntheorem}

\newtheorem{mdef}{Definition}[section]
\newtheorem{rmk}{Remark}

\newtheorem{thm}{Theorem}[section]
\newtheorem{thmm}{Theorem}
\newtheorem{prop}[thm]{Propositon}
\newtheorem{lem}[thm]{Lemma}
\newtheorem{cor}[thm]{Corollary}
\newtheorem{Cor}{Corollary}[thmm]
\newtheorem{fact}{Fact}

{
\theoremstyle{nonumberchange}
\theoremsymbol{\mbox{$\Box$}}
\newtheorem{pf}{Proof}
\newtheorem{potA}{Proof of Theorem \ref{thm:localvsfinite}}
\newtheorem{potB}{Proof of Theorem \ref{thm:SMB}}
\newtheorem{potC}{Proof of Theorem \ref{thm:vp}}
\newtheorem{pof1}{Proof of Fact \ref{fact:convergence}}
}

\title{Unstable Entropy and Unstable Pressure for Partially Hyperbolic Endomorphisms\footnotetext{
X.Wang is the corresponding author.}}
\author{Xinsheng Wang$^1$, Weisheng Wu$^2$ and Yujun Zhu$^3$}
\date{}
\begin{document}
\maketitle
{\footnotesize
  \centerline{1. School of Mathematical Sciences}
  \centerline{Hebei Normal University, Shijiazhuang, 050024, P.R. China}
  \centerline{2. Department of Applied Mathematics, College of Science}
  \centerline{China Agricultural University, Beijing, 100083, P.R. China}
  \centerline{3. School of Mathematical Sciences}
  \centerline{Xiamen University, Xiamen, 361005, P.R. China}
}

\begin{abstract}

In this paper, unstable metric entropy, unstable topological entropy and unstable pressure for partially hyperbolic endomorphisms are introduced and investigated. A version of Shannon-McMillan-Breiman Theorem is established, and a variational principle is formulated, which gives a relationship between unstable metric entropy and unstable pressure (unstable topological entropy). As an application of the variational principle, some results on the $u$-equilibrium states are given.

\end{abstract}

\section{Introduction}
In order to describe the complexity of a dynamical system from different points of view, some invariants are introduced, among which, entropy including metric entropy and topological entropy is a crucial one. The metric entropy gives the maximum amount of average information one can get from a system with respect to an invariant measure, while the topological entropy describes the exponential growth rate of the number of orbits. A variational principle relating them says that topological entropy is equal to  the supremum of metric entropy over all invariant measures. The pressure with respect to a potential function is a generalization of topological entropy, and a variational principle relating metric entropy and it can also be formulated.

Concerning smooth ergodic theory, the Lyapunov exponent can be introduced and there has been a system of results related with entropy for diffeomorphisms (the invertible case). Pesin's entropy formula relating metric entropy and Lyapunov exponents with respect to an SRB measure is established for both deterministic and random cases \textcolor[rgb]{0.00,0.00,0.00}{(cf.  \cite{LedrappierYoung1985},  \cite{LiuQian1995})}. Moreover, a generalized entropy formula (dimension formula) with respect to a general invariant measure is formulated (cf.  \cite{LedrappierYoung1985a}).

Plenty of physical processes are irreversible, so it is interesting to investigate corresponding results as in \cite{LedrappierYoung1985,LedrappierYoung1985a} for non-invertible case, i.e. endomorphisms. However, for endomorphisms, there are some difficulties to establish similar results.  Due to non-invertibility, the preimage of a given point is usually not a single point, hence the notion of unstable manifolds is not well defined and therefore it is a problem to formulate the SRB property. Further more, the non-invertibility leads to some subtle difficulties such that it is not convenient to consider problems concerning entropies on the phase space. In order to overcome this difficulty, in \cite{Zhu1998}, Zhu introduced the inverse limit space (see Section \ref{sec:pre}, for details), which makes it possible to define the unstable manifolds and borrow some ideas from the smooth ergodic theory for random dynamical systems. In \cite{QianZhu2001}, Qian and Zhu gave the necessary and sufficient condition for Pesin's entropy formula in the case of endomorphisms. Moreover, a series of results on ergodic theory of endomorphisms are obtained, see  \cite{QianXieZhu2009} for a complete discussion of this topic.

Via the inverse limit space, the unstable manifolds can be well defined, which implies us that the unstable entropy and unstable pressure can be introduced for endomorphisms. The concept of unstable entropy is originally introduced by Hu, Hua and Wu in  \cite{HuHuaWu2017} for partially hyperbolic diffeomorphisms, which gives a description of the complexity of a system along unstable manifolds. In  \cite{HuHuaWu2017}, a complete discussion is given, including the relationship between unstable metric entropy and Ledrappier and Young's entropy, a version of Shannon-McMillan-Breiman Theorem and a variational principle relating unstable metric entropy and unstable topological entropy. It is important to point out that in  \cite{HuHuaWu2017} the unstable metric entropy is given by the conditional entropy of a finite partition with respect to a measurable partition, instead of the form in Ledrappier and Young's papers \cite{LedrappierYoung1985,LedrappierYoung1985a}. The latter form of metric entropy is introduced to establish connection with SRB measures, which is not suitable for giving an unstable version of variational principle; while the former one makes it possible to formulate the variational principle using a classical method. However, both of the two forms give the same thing but from different points of view (see Theorem A in  \cite{HuHuaWu2017} for more details). As a generalization of unstable topological entropy, unstable pressure is defined in  \cite{HuWuZhu2017}, where the so-called $u$-equilibrium states are introduced and investigated finely. Recently, a version of local variational principle and Katok's entropy formula for unstable metric entropy are given in \cite{Wu2017} and \cite{HuangChenWang2018} respectively.

Our main purpose in this paper is to establish unstable entropy and unstable pressure for endomorphisms and try to give a system of results as in \cite{HuHuaWu2017}. Inspired by the argument in  \cite{QianZhu2001,QianXieZhu2009}, we introduce the concepts via the inverse limit space.

For an endomorphism $f$ on a Riemannian manifold $M$ with an $f$-invariant measure $\mu$, the inverse limit space $M^f\subset M^{\mathbb{Z}}$ is introduced, and a dynamical system $(M^f,\tau,\tilde{\mu})$ is established, where $\tau$ is the left shift operator on $M^f$ and $\tilde{\mu}$ is a $\tau$-invariant measure corresponding to $\mu$. Thanks to $(M^f,\tau,\tilde{\mu})$, we can give two types of definitions of unstable metric entropy, one is introduced via a ``pointwise'' way (see Definition \ref{def:metricentropy1}), which is denoted by $\tilde{h}^u_\mu(f)$, and the other one is defined using finite partitions (see Definition \ref{def:metricentropy2}), which is in the classical form and denoted by $h^u_\mu(f)$. Then we show that $\tilde{h}^u_\mu(f)$ and $h^u_\mu(f)$ are equivalent when $\tilde{\mu}$ is ergodic (Theorem \ref{thm:localvsfinite}). Then a version of Shannon-McMillan-Breiman Theorem is established for our case (Theorem \ref{thm:SMB}), which makes our unstable metric entropy meaningful. Again using $(M^f,\tau,\tilde{\mu})$, we define the unstable pressure and unstable topological entropy, the latter one can be viewed as a special case of the former one. And we show that the variational principle for classical entropy and pressure also holds in our case (Theorem \ref{thm:vp}), which makes it possible to consider the so-called $u$-equilibrium states for endomorphisms.

This paper is organized as follows. In Section \ref{sec:pre}, we give some preliminaries, including the concept of partial hyperbolicity for endomorphisms, the inverse limit space and other necessary definitions for our results. And in the end of this section, our main results are also given. In Section \ref{sec:umetricentropy}, we give the precise definitions of unstable metric entropy for endomorphisms via two methods, following which, Theorem \ref{thm:localvsfinite} is proved. And Theorem \ref{thm:SMB} is also proved in this section. In section \ref{sec:pressure}, definitions of unstable pressure and unstable topological entropy are introduced, some properties of unstable pressure are also listed in the end of this section. In section \ref{sec:vp}, we prove our main result, i.e., the variational principle (Theorem \ref{thm:vp}), and as an application, the so-called $u$-equilibrium states for endomorphisms are introduced.
\section{Preliminaries and statements of main results}\label{sec:pre}
Throughout this paper, let $M$ be a compact $C^\infty$ Riemannian manifold without boundary endowed with metric $d(\cdot,\cdot)$ and $f:M\to M$ a $C^1$ endomorphism. Denote $TM$ the tangent bundle of $M$ with norm $\Vert\cdot\Vert$. Both $d(\cdot,\cdot)$ and $\Vert\cdot\Vert$ are induced by the Riemannian metric.

  For a metric space $X$, denote $\mathcal{B}(X)$ the Borel $\sigma$-algebra of $X$. Let $M^\mathbb{Z}$ be the infinite product space of $M$ endowed with the product topology and the metric $\tilde{d}(\tilde{x},\tilde{y})=\sum_{n=-\infty}^{\infty}2^{-|n|}d(x_n,y_n)$ for $\tilde{x}=\{x_n\}_{n=-\infty}^{\infty}$ and $\tilde{y}=\{y_n\}_{n=-\infty}^{\infty}$. In order to define unstable manifolds, we need the concept of \textit{inverse limit space} denoted by $M^f$, which means it is a subspace of the product space $M^\mathbb{Z}$, and $fx_n=x_{n+1}$, $n\in\mathbb{Z}$, for $\tilde{x}=\{x_n\}_{n=-\infty}^{+\infty}\in M^f$. It is clear that $M^f$ is a closed subspace of $M^{\mathbb{Z}}$. Let $\Pi\colon M^f\to M$ be the projection such that for $\tilde{x}=\{x_n\}_{n=-\infty}^{+\infty}$, $\Pi(\tilde{x})=x_0$.

  Let $\tau\colon M^f\to M^f$ be the left shift operator. Denote $\mathcal{M}(f)$ the set of all $f$-invariant Borel measures on $M$, and denote $\mathcal{M}(\tau)$ the set of all $\tau$-invariant measures on $M^f$. On one hand, for any $\mu\in\mathcal{M}(f)$, there is a unique $\tau$-invariant measure $\tilde{\mu}$ on $M^f$ corresponding to $\mu$ with $\Pi(\tilde{\mu})=\mu$ (see Proposition \uppercase\expandafter{\romannumeral1}.3.1 in  \cite{QianXieZhu2009}); on the other hand, for any $\tilde{\mu}\in\mathcal{M}(\tau)$, $\mu\colon=\Pi(\tilde{\mu})$ is an $f$-invariant measure on $M$. In addition, $\mu$ is ergodic with respect to $f$ if and only if $\tilde{\mu}$ is ergodic with respect to $\tau$. For more details on the relationship between $\mathcal{M}(f)$ and $\mathcal{M}(\tau)$, the reader can refer to \uppercase\expandafter{\romannumeral1}.3 in  \cite{QianXieZhu2009}. In the remaining of this paper, we always denote $\mu$ and $\tilde{\mu}$ the measures on $M$ and $M^f$ respectively with $\Pi(\tilde{\mu})=\mu$.

  Consider the pull back bundle $E=\Pi^*TM$. The tangent map $Df$ induces a fiber preserving map on $E$ with respect to the left shift map $\tau$, defined by $\Pi^*\circ Df\circ\Pi_*$, and still denoted by $Df$ for simplicity. Now, we give the definition of partial hyperbolicity.
\begin{mdef}\rm\label{def:ph}
  \textcolor[rgb]{0.00,0.00,0.00}{$f$ is said to be \textit{(uniformly) partially hyperbolic} if there exist a continuous splitting of the pull back bundle $E$ into three subbundles, i.e., $E(\tilde{x})=E^s({\tilde{x}})\oplus E^c({\tilde{x}})\oplus E^u({\tilde{x}})$ for all ${\tilde{x}}\in M^f$ and constants $\lambda_1$, $\lambda_1'$, ${\lambda_2}$, ${\lambda_2}'$ and $C$ with $0<{\lambda_1}<1<{\lambda_2}$, ${\lambda_1}<{\lambda_1}'\leq{\lambda_2}'<{\lambda_2}$ and $C>0$ such that for each $\tilde{x}\in M^f$,}
  \begin{enumerate}[label=(\roman*)]
    \item {$D_{\tilde{x}}f(E^k({\tilde{x}}))= E^k(\tau({\tilde{x}}))$}, for $k=s,c,u$;
    \item for $v^s\in E^s({\tilde{x}})$ and $n\in\mathbb{Z}^+$, $\Vert D_{\tilde{x}}f^nv^s\Vert\leq C{\lambda_1}^n\Vert v^s\Vert$;
    \item for $v^c\in E^c({\tilde{x}})$ and $n\in\mathbb{Z}^+$, $C^{-1}{({\lambda_1}')}^{n}\Vert v^c\Vert\leq\Vert D_{\tilde{x}}f^nv^c\Vert\leq C({\lambda_2}')^n\Vert v^c\Vert$;
    \item for $v^u\in E^u({\tilde{x}})$ and $n\in\mathbb{Z}^+$, $\Vert D_{\tilde{x}}f^nv^u\Vert\geq C^{-1}{{\lambda_2}}^{n}\Vert v^u\Vert$.
  \end{enumerate}
\end{mdef}

  From now on, let $(f,M,\mu)$ be a dynamical system, where $f$ is a partially hyperbolic endomorphism, and $\mu$ is an $f$-invariant Borel measure. Let $\tilde{\mu}$ be the corresponding measure on $M^f$.

  For $\tilde{x}=\{x_n\}_{n=-\infty}^{\infty}\in M^f$ and $\epsilon>0$ small enough, define
  \begin{align*}
    W^u_\epsilon(\tilde{x},f)\colon=\{&z_0\in M\colon \text{there exists }\tilde{z}\in M^f \text{ with }\Pi(\tilde{z})=z_0,\\&d(z_{-n},x_{-n})<\epsilon\text{ for } n\in\mathbb{N} \text{ and }\limsup_{n\to\infty}\frac{1}{n}\log d(z_{-n},x_{-n})\leq-\log\lambda_2\},
  \end{align*}
  where $\lambda_2$ is the constant in Definition \ref{def:ph}. $W^u_\epsilon(\tilde{x},f)$ is called a \textit{local unstable manifold} of $f$ at $\tilde{x}$. Now we have the following theorem, which is stated for hyperbolic endomorphisms, while it is still valid for our partially hyperbolic case. The reader can also refer to \cite{Przytycki1976,Ruelle1988,Young1986} for more details.

  \begin{thm}[Theorem \uppercase\expandafter{\romannumeral4}.2.1 in \cite{QianXieZhu2009}]\label{thm:umt}
    Let $f$ be a partially hyperbolic endomorphism. Then there exists a continuous family of $C^1$ embedded  disks $\{D^u_{\tilde{x}}\}_{\tilde{x}\in M^f}$ in $M$ and constants $0<\lambda<1$ and $\epsilon>0$ such that
    \begin{enumerate}[label=(\roman*)]
      \item $T_{{x}_0}D^u_{\tilde{x}}=E^u({x}_0)$, for any $\tilde{x}\in M^f$;
      \item for any $z_0\in D^u_{\tilde{x}}$, there exists unique $\tilde{z}\in M^f$ such that $\Pi(\tilde{z})=z_0$ and
          \begin{equation}\label{eq:locunstable2}
            d(z_{-n},x_{-n})\leq \lambda^nd(z_0,x_0),
          \end{equation}
      for $n\in\mathbb{Z}^+$;
      \item $D^u_{\tilde{x}}\cap B(x_0,\epsilon)=W^u_\epsilon(\tilde{x},f)$, where $B(x_0,\epsilon)=\{y\in M\colon d(y,x)<\epsilon\}$.
    \end{enumerate}
  \end{thm}
  Then we can define
  \[
    \widetilde{W}^u_\epsilon(\tilde{x},f)\colon=\{\tilde{z}\in M^f\colon\Pi(\tilde{z})\in W^u_\epsilon(\tilde{x},f)\text{ and }\tilde{z}\text{ satisfies \eqref{eq:locunstable2}}\}.
  \]
  Sometimes, we will use the notation $W^u_{\text{loc}}(\tilde{x},f)$ and $\widetilde{W}^u_{\text{loc}}(\tilde{x},f)$ for $W^u_{\epsilon}(\tilde{x},f)$ and $\widetilde{W}^u_{\epsilon}(\tilde{x},f)$ respectively.
  \begin{rmk}\rm
    According to Theorem \ref{thm:umt}, it is clear that
    \[
      \Pi|_{\widetilde{W}^u_{\text{loc}}(\tilde{x},f)}\colon\widetilde{W}^u_{\text{loc}}(\tilde{x},f)\to W^u_{\text{loc}}(\tilde{x},f)
    \]
    is a bijection, which is crucial for our subsequent proofs.
  \end{rmk}
  Now we define
  \begin{align*}
    W^u(\tilde{x},f)=\{&z_0\in M\colon \text{there exits }\tilde{z}\text{ with }\Pi(\tilde{z})=z_0\\
    &\text{and }\textcolor[rgb]{0.00,0.00,0.00}{\limsup_{n\to+\infty}\frac{1}{n}d(z_{-n},x_{-n})\leq-\log\lambda_2}\},
  \end{align*}
  and
  \[
    \widetilde{W}^u(\tilde{x},f)=\{\tilde{z}\in M^f\colon\Pi(\tilde{z})\in W^u(\tilde{x},f)\text{ with }{\limsup_{n\to+\infty}\frac{1}{n}d(z_{-n},x_{-n})\leq-\log\lambda_2}\},
  \]
  where $\lambda_2$ is the constant in Definition \ref{def:ph}. We call $W^u(\tilde{x},f)$ the \textit{global unstable set} at $\tilde{x}$.

  It can also be proved as in  \cite{Zhu1998} that there exists a sequence of $C^1$ embedded disks $\{W_{-n}(\tilde{x})\}_{n=0}^{+\infty}$ in $M$ such that $fW_{-n}(\tilde{x})\supset W_{-(n-1)}(\tilde{x})$ for $n\in\mathbb{Z}^+$ and
  \[
    W^u(\tilde{x},f)=\bigcup_{n=0}^{+\infty}f^nW_{-n}(\tilde{x}),
  \]
  which shows that $W^u(\tilde{x},f)$ is in fact an immersed submanifold of $M$ tangent at $\Pi(\tilde{x})$ to $E^u(\Pi(\tilde{x}))$. Then we denote the set $\{W^u(\tilde{x},f)\colon\tilde{x}\in M^f\}$ by $W^u$, which is called $W^u$-foliation.

  For a measurable partition $\eta$ of $M^f$, $\eta(\tilde{x})$ means the element in $\eta$ containing $\tilde{x}$. Now we give some definitions related to measurable partitions.

  \begin{mdef}\rm
    A measurable partition $\eta$ of $M^f$ is said to \textit{be subordinate to} $W^u$-foliation if for $\tilde{\mu}$-a.e. $\tilde{x}$, $\eta(\tilde{x})$ has the following properties:
    \begin{enumerate}[label=(\roman*)]
      \item $\Pi|_{\eta{(\tilde{x})}}\colon\eta(\tilde{x})\to\Pi(\eta(\tilde{x}))$ is bijective;
      \item There exists a $k(\tilde{x})$-dimensional \textcolor[rgb]{0.00,0.00,0.00}{(where $k(\tilde{x})=\dim  E^u(x_0)$)} $C^1$ embedded submanifold $W_{\tilde{x}}$ of $M$ with $W_{\tilde{x}}\subset W^u(\tilde{x})$, such that $\Pi(\eta(\tilde{x}))\subset W_{\tilde{x}}$, and $\Pi(\eta(\tilde{x}))$ contains an open neighborhood of $x_0$ in $W_{\tilde{x}}$.
    \end{enumerate}
  \end{mdef}

Given $\tilde{\mu}\in\mathcal{M}(\tau)$. For a measurable partition $\eta$ of $M^f$, there exists a canonical system $\{\tilde{\mu}^\eta_{\tilde{x}}\}_{\tilde{x}\in M^f}$ of conditional measures of $\tilde{\mu}$ associated with $\eta$, satisfying
\begin{enumerate}[label=(\roman*)]
  \item for every measurable set $\widetilde{B}\subset M^f$, $\tilde{x}\mapsto\tilde{\mu}^\eta_{\tilde{x}}(\widetilde{B})$ is measurable;
  \item $\tilde{\mu}(\widetilde{B})=\int_{M^f}\tilde{\mu}^\eta_{\tilde{x}}(\widetilde{B})\mathrm{d}\tilde{\mu}(\tilde{x})$.
\end{enumerate}
See e.g.  \cite{Rokhlin1962} for more details.

Let $\alpha$ be a measurable partition of $M^f$. The diameter of $\alpha$ is defined as follows:
\[
  \mathrm{diam}(\alpha)=\sup_{A\in\alpha}\mathrm{diam}(\Pi(A)),
\]
where for a subset $B$ of $M$,
\[
  \mathrm{diam}(B):=\sup_{x,y\in B}d(x,y).
\]

In the following, we construct a type of measurable partitions subordinate to $W^u$-foliation. Fix $\epsilon>0$. Let $\alpha$ be a finite partition of $M^f$ with diameter small enough such that the diameter of $\alpha$ is less than $\epsilon$. We can construct a finer partition $\eta$ such that for each $\tilde{x}\in M^f$
\[
  \eta(\tilde{x})=\alpha(\tilde{x})\cap\widetilde{W}^u_\epsilon(\tilde{x},f).
\]
Clearly $\eta$ is a measurable partition of $M^f$ (cf. p34 in \cite{HuHuaWu2017} for more details). In addition, by the definition of $\widetilde{W}^u_\epsilon(\tilde{x},f)$ and Theorem \ref{thm:umt}, if $\mu(\partial(\Pi(\alpha)))=0$, $\eta$ is a measurable partition subordinate to $W^u$-foliation, where $\partial(\Pi(\alpha))=\bigcup_{A\in\alpha}\partial(\Pi (A))$ and for $B\subset M$, $\partial B$ means the boundary of $B$. Let $\mathcal{P}(M^f)$ denote the set of all finite partitions with diameter small enough and $\mathcal{P}^u(M^f)$ denote the set of measurable partitions of $M^f$ subordinate to $W^u$-foliation which are induced by finite partitions in $\mathcal{P}(M^f)$.

In the following, we consider a special type of measurable partitions.
\begin{mdef}\rm
  A measurable partition $\xi$ of $M^f$ is said to be \textit{increasing} if $\tau^{-1}\xi\geq\xi$.
\end{mdef}

Consider a measurable partition $\xi=\{A_i\}_{i\in I}$ of $M^f$. A measurable set $B$ is called a $\xi$-set if $B=\cup_{i\in I'}A_i$, where $I'\subset I$. Denote $\mathcal{B}(\xi)$ the $\sigma$-algebra of $M^f$ consisting of all measurable $\xi$-sets. Given $\tilde{\mu}\in\mathcal{M}(\tau)$, define
\[
  \mathcal{B}^u\colon=\{\tilde{B}\in\mathcal{B}_{\tilde{\mu}}(M^f)\colon\tilde{x}\in \tilde{B}\text{ implies }\widetilde{W}^u(\tilde{x})\subset \tilde{B}\},
\]
where $\mathcal{B}_{\tilde{\mu}}(M^f)$ is the completion of $\mathcal{B}(M^f)$ with respect to $\tilde{\mu}$.

The following proposition ensures the existence of increasing measurable partitions, the reader can see Section {\uppercase\expandafter{\romannumeral9}}.2.2 in \cite{QianXieZhu2009} for details.

\begin{prop}\label{prop:specialpartition}
  There exists a measurable partition $\xi$ of $M^f$ which has the following properties:
  \begin{enumerate}[label=(\roman*)]
    \item $\tau^{-1}\xi\geq\xi$;
    \item $\bigvee^\infty_{n=0}\tau^{-n}\xi$ is equal to the partition into single points;
    \item $\mathcal{B}(\bigwedge^\infty_{n=0}\tau^n(\xi))=\mathcal{B}^u$, $\tilde{\mu}$-$\mathrm{mod}$ $0$;
    \item \textcolor[rgb]{0.00,0.00,0.00}{$\xi$ is subordinate to $W^u$-foliation of $f$.}
  \end{enumerate}
\end{prop}

We denote by $\mathcal{Q}^u(M^f)$ the set of all increasing measurable partitions subordinate to $W^u$-foliation as in Proposition \ref{prop:specialpartition}.

Now we can introduce the unstable metric entropy along $W^u$-foliation for a partially hyperbolic endomorphism. Two types of unstable metric entropy will be given, one is defined via the average decreasing rate of the Bowen balls which is denoted by $\tilde{h}^u_\mu(f)$, and the other one is defined by the conditional entropy of $f$ along $W^u$-foliation, which is denoted by $h^u_\mu(f)$. Both of their precise definitions are stated in Section \ref{sec:umetricentropy}. When $\mu$ is ergodic, it can be showed that $\tilde{h}^u_\mu(f)$ and $h^u_\mu(f)$ describe the same thing from different points of view essentially, i.e., we have the following theorem.
\begin{thmm}\label{thm:localvsfinite}
Let $f$ be a $C^1$ partially hyperbolic endomorphism and $\mu$ an ergodic measure of ${f}$. Then
  \[
    \tilde{h}^u_\mu({f})={h}^u_\mu({f}).
  \]
\end{thmm}

The classical Shannon-McMillan-Breiman Theorem expresses the metric entropy as the limit of certain conditional information functions, which interprets the metric entropy from the viewpoint of information theory. Moreover, there is a corresponding version of Shannon-McMillan-Breiman Theorem in our case. Specially, for two measurable partitions $\beta$ and $\gamma$ of $M^f$, the conditional entropy of $\beta$ with respect to $\gamma$ for a $\tau$-invariant measure $\tilde{\mu}$ can be given, which is denoted by $I_{\tilde{\mu}}(\beta|\gamma)$. Meanwhile, the conditional entropy $h_\mu(f,\beta|\gamma)$ of $f$ for $\beta$ with respect to $\gamma$ can be introduced. Precise definitions of above concepts are stated in Section \ref{sec:umetricentropy}. Then we have the following theorem.
\begin{thmm}\label{thm:SMB}
 Let $f$ be a $C^1$ partially hyperbolic endomorphism, and suppose $\mu$ is an ergodic measure of ${f}$. {For any $\alpha\in \mathcal{P}( M^f), \eta\in\mathcal{P}^u( M^f)$}, we have
  \[
    \lim_{n\to\infty}\frac{1}{n}I_{\tilde{\mu}}(\alpha^{n-1}_0|\eta)(\tilde{x})=h_\mu({f},\alpha|\eta),
  \]
  where for integers $k<j$, $\tau^{-j}\beta\vee\tau^{-(j-1)}\beta\vee\cdots\vee\tau^{-k}\beta$ is denoted by $\beta^j_k$.
\end{thmm}

Taking account of the topological structure, we can establish the concepts of unstable topological entropy and unstable pressure with respect to a continuous potential function $\varphi$, whose definitions will be given in Section \ref{sec:pressure}. We denote unstable topological entropy and unstable pressure by $h^u_{\text{{\rm top}}}(f)$ and $P^u(f,\varphi)$ respectively.

It is natural to consider the relationship between unstable pressure and unstable metric entropy, a version of variational principle can be formulated as follows. Denote $C(M)$ the set of all continuous functions on $M$.

\begin{thmm}\label{thm:vp}
  Let $f$ be a $C^1$ partially hyperbolic endomorphism and $\varphi\in C(M)$. Then we have
  \[
    \sup_{\mu\in\mathcal{M}({f})}\left\{h^u_\mu({f})+\int_{ M}\varphi\mathrm{d}{\mu}\right\}=P^u({f},\varphi).
  \]
\end{thmm}

A direct corollary of Theorem \ref{thm:vp} is the following variational principle for unstable topological entropy.
\begin{Cor}\label{coro:vp1}
Let $f$ be a $C^1$ partially hyperbolic endomorphism, then we have
 \[
   \sup_{\mu\in\mathcal{M}({f})}\left\{h^u_\mu({f})\right\}=h_{\text{{\rm top}}}^u({f}).
 \]
\end{Cor}

\section{Unstable metric entropy for endomorphisms}\label{sec:umetricentropy}
\subsection{Definitions of unstable metric entropy}

In this subsection, we give the definition of unstable metric entropy for endomorphisms via two methods. The equivalence of the two definitions will be proved in the next subsection. Firstly, we give the definition by ``pointwise'' approach.

Let $d^u_{\tilde{x}}$ be the metric on $W^u(\tilde{x},f)$ induced by the Riemannian structure on $M$. Denote the $d^u_{n}$-Bowen ball in $\widetilde{W}^u(\tilde{x},f)$ with center $\tilde{x}$ and radius $\epsilon>0$ by $\widetilde{V}^u(f,\tilde{x},n,\epsilon)$, i.e.,
\[
  \widetilde{V}^u(f,\tilde{x},n,\epsilon):=\{\tilde{y}\in \widetilde{W}^u(\tilde{x},f)\colon \tilde{d}^u_{n}(\tilde{y},\tilde{x})<\epsilon\},
\]
where
\[
   \tilde{d}^u_{n}(\tilde{x},\tilde{y}):=\max_{0\leq j\leq n-1}\{ d^{u}_{\tau^j\tilde{x}}(\Pi(\tau^j\tilde{x}),\Pi(\tau^j\tilde{y}))\}.
\]
\begin{mdef}\rm\label{def:metricentropy1}
  Given an increasing partition $\xi_u$ of $M^f$ subordinate to $W^u$-foliation, we define the \textit{unstable metric entropy} along $W^u$-foliation as follows:
\[
  h_\mu(f,\xi_u)=\int_{M^f}h_\mu(f,\tilde{x},\xi_u)\mathrm{d}{\tilde\mu}(\tilde{x}),
\]
where
\[
  h_\mu(f,\tilde{x},\xi_u)=\lim_{\epsilon\to 0}\limsup_{n\to\infty}-\frac{1}{n}\log{\tilde{\mu}}_{\tilde{x}}^{\xi_u}\widetilde{V}^u(f,\tilde{x},n,\epsilon).
\]
\end{mdef}

It can be proved that $h_\mu(f,\tilde{x},\xi_u)$ is independent of the choice of $\xi_u$, hence we also denote $\tilde{h}^u_\mu(f)= h_\mu(f,\xi_u)$. Moreover, $h_\mu(f,\tilde{x},\xi_u)$ is $\tau$-invariant, so when $\tilde{\mu}$ is ergodic, we have $\tilde{h}^u_\mu(f)=h_\mu(f,\xi_u)=h_\mu(f,\tilde{x},\xi_u)$, for $\tilde{\mu}$-a.e. $\tilde{x}\in M^f$. The reader can refer to Section {\uppercase\expandafter{\romannumeral9}}.3 in \cite{QianXieZhu2009} for details.

\begin{rmk}\rm\label{rmk:low=up}
  In fact, in Definition \ref{def:metricentropy1}, ``$\limsup$'' can be replaced by ``$\lim$'' and ``$\lim_{\epsilon\to 0}$'' can be dropped. Indeed, denote
  \begin{equation*}\label{eq:lower}
    \underline{h}({f},\tilde{x},\epsilon,\xi_u)=\liminf_{n\to\infty}-\frac{1}{n}\log{\tilde\mu}_{\tilde{x}}^{\xi_u}\widetilde{V}^u({f},\tilde{x},n,\epsilon)
  \end{equation*}
  and
  \begin{equation*}\label{eq:upper}
    \overline{h}({f},\tilde{x},\epsilon,\xi_u)=\limsup_{n\to\infty}-\frac{1}{n}\log{\tilde\mu}_{\tilde{x}}^{\xi_u}\widetilde{V}^u({f},\tilde{x},n,\epsilon).
  \end{equation*}
  It has been proved in Section {\uppercase\expandafter{\romannumeral9}}.3 of \cite{QianXieZhu2009} that
  \begin{equation*}\label{eq:low=up}
    \lim_{\epsilon\to0}\underline{h}({f},\tilde{x},\epsilon,\xi_u)=\lim_{\epsilon\to0}\overline{h}({f},\tilde{x},\epsilon,\xi_u).
  \end{equation*}
  Then following the proof of Lemma 3.1 in  \cite{HuHuaWu2017}, we can prove the above claim since $f$ is uniformly expanding along $W^u$-foliation.
\end{rmk}

In order to give the definition of unstable metric entropy via measurable partitions, firstly we give some definitions on information function, which is slightly modified in our context. Some properties concerning information function will be also listed in the end of this subsection.

\begin{mdef}\rm
  Let $\alpha$ and $\eta$ be two measurable partitions of $M^f$. The \textit{information function} of $\alpha$ with respect to ${\tilde\mu}$ is defined as
  \[
    I_{{\tilde{\mu}}}(\alpha)(\tilde{x}):=-\log{\tilde{\mu}}(\alpha(\tilde{x})),
  \]
  and the \textit{entropy} of $\alpha$ with respect to ${\tilde\mu}$ is defined as
  \[
    H_{{\tilde{\mu}}}(\alpha):=\int_{M^f}I_{{\tilde{\mu}}}(\alpha)(\tilde{x})\mathrm{d}{\tilde{\mu}}(\tilde{x})=-\int_{M^f}\log{\tilde{\mu}}(\alpha(\tilde{x}))\mathrm{d}{\tilde{\mu}}(\tilde{x}).
  \]
  The \textit{conditional information function} of $\alpha$ with respect to $\eta$ is defined as
  \[
    I_{{\tilde{\mu}}}(\alpha|\eta)(\tilde{x}):=-\log{\tilde{\mu}}_{\tilde{x}}^\eta(\alpha(\tilde{x})),
  \]
 where $\{{\tilde{\mu}}^\eta_{(\tilde{x})}\}_{\tilde{x}\in M^f}$ is a canonical system of conditional measures of ${\tilde{\mu}}$ with respect to $\eta$. Then the \textit{conditional entropy} of $\alpha$ with respect to $\eta$ is defined as
   \[
    H_{{\tilde{\mu}}}(\alpha|\eta):=\int_{M^f}I_{{\tilde{\mu}}}(\alpha|\eta)(\tilde{x})\mathrm{d}{\tilde{\mu}}(\tilde{x})=-\int_{M^f}\log{\tilde{\mu}}_{\tilde{x}}^\eta(\alpha(\tilde{x}))\mathrm{d}{\tilde{\mu}}(\tilde{x}).
  \]
\end{mdef}

Now we can give the definition of unstable metric entropy by finite partitions.

\begin{mdef}\rm\label{def:metricentropy2}
  The \textit{conditional entropy} of $f$ for a finite measurable partition $\alpha$ of $M^f$ with respect to $\eta\in\mathcal{P}^u(M^f)$ is defined as
  \[
    h_\mu(f,\alpha|\eta)=\limsup_{n\to\infty}\frac{1}{n}H_{\tilde{\mu}}(\alpha^{n-1}_0|\eta).
  \]
  \textit{The conditional entropy} of $f$ with respect to $\eta$ is defined as
  \[
    h_{\mu}(f|\eta)=\sup_{\alpha\in\mathcal{P}(M^f)}h_\mu(f,\alpha|\eta),
  \]
  and the \textit{conditional entropy} of $f$ along $W^u$-foliation is defined as
  \[
    {h}^u_\mu(f)=\sup_{\eta\in\mathcal{P}^u(M^f)}h_{\mu}(f|\eta).
  \]
\end{mdef}

To end this subsection, we list the following lemmas which are derived from  \cite{HuHuaWu2017} with slight adaption and will be useful for the proofs of our main results.
\begin{lem}\label{lem:info}
  Given $\tilde\mu\in\mathcal{M}({\tau})$ and let $\alpha$, $\beta$ and $\gamma$ be measurable partitions of $ M^f$ with $H_{\tilde{\mu}}(\alpha|\gamma)$, $H_{\tilde{\mu}}(\beta|\gamma)<\infty$.
  \begin{enumerate}[label=(\roman*)]
    \item If $\alpha\leq\beta$, then $I_{\tilde{\mu}}(\alpha|\gamma)(\tilde{x})\leq I_{\tilde{\mu}}(\beta|\gamma)(\tilde{x})$ and $H_{\tilde{\mu}}(\alpha|\gamma)\leq H_{\tilde{\mu}}(\beta|\gamma)$;
    \item$I_{\tilde{\mu}}(\alpha\vee\beta|\gamma)(\tilde{x})=I_{\tilde{\mu}}(\alpha|\gamma)(\tilde{x})+I_{\tilde{\mu}}(\beta|\alpha\vee\gamma)(\tilde{x})
    $
    and $H_{\tilde{\mu}}(\alpha\vee\beta|\gamma)=H_{\tilde{\mu}}(\alpha|\gamma)+H_{\tilde{\mu}}(\beta|\alpha\vee\gamma)$;
    \item $H_{\tilde{\mu}}(\alpha\vee\beta|\gamma)\leq H_{\tilde{\mu}}(\alpha|\gamma)+H_{\tilde{\mu}}(\beta|\gamma);$
    \item if $\beta\leq\gamma$, then $H_{\tilde{\mu}}(\alpha|\beta)\geq H_{\tilde{\mu}}(\alpha|\gamma)$.
  \end{enumerate}
\end{lem}
\begin{lem}\label{lem:info2}
Let $\tilde\mu\in\mathcal{M}({\tau})$, and $\alpha$, $\beta$ and $\gamma$ measurable partitions of $ M^f$.
  \begin{enumerate}[label=(\roman*)]
    \item\label{lem:info2_item1}
    \[
      I_{\tilde{\mu}}(\beta_0^{n-1}|\gamma)(\tilde{x})=I_{\tilde{\mu}}(\beta|\gamma)(\tilde{x})+\sum_{i=1}^{n-1}I_{\tilde{\mu}}(\beta|\tau^i(\beta_0^{i-1}\vee\gamma))(\tau^i(\tilde{x})),
    \]
    hence
    \[
      H_{\tilde{\mu}}(\beta_0^{n-1}|\gamma)=H_{\tilde{\mu}}(\beta|\gamma)+\sum_{i=1}^{n-1}H_{\tilde{\mu}}(\beta|\tau^i(\beta_0^{i-1}\vee\gamma));
    \]
    \item\label{lem:info2_item2}
    \begin{align*}
       & I_{\tilde{\mu}}(\alpha_0^{n-1}|\gamma)(\tilde{x}) \\
      =& I_{\tilde{\mu}}(\alpha|\tau^{n-1}\gamma)(\tau^{n-1}(\tilde{x}))+\sum_{i=0}^{n-2}I_{\tilde{\mu}}(\alpha|\alpha_1^{n-1-i}\vee \tau^i\gamma)(\tau^i(\tilde{x})),
    \end{align*}
    hence
    \[
      H_{\tilde{\mu}}(\alpha_0^{n-1}|\gamma)=H_{\tilde{\mu}}(\alpha|\tau^{n-1}\gamma)+\sum_{i=0}^{n-2}H_{\tilde{\mu}}(\alpha|\alpha_1^{n-1-i}\vee \tau^i\gamma).
    \]
  \end{enumerate}
\end{lem}
\begin{lem}\label{lem:semicontiofparti}
  Let $\alpha\in\mathcal{P}( M^f)$ and $\{\zeta_n\}$ be a sequence of increasing measurable partitions of $ M^f$ with $\zeta_n\nearrow\zeta$. Then for $\varphi_n(\tilde{x})=I_{\tilde{\mu}}(\alpha|\zeta_n)(\tilde{x})$, $\varphi^*:=\sup_n\varphi_n\in L^1({\mu})$.
   %where $L^1({\mu})$  is the set of all $L^1$ functions on $ M^f$ with respect to ${\mu}$.
\end{lem}
\begin{lem}\label{lem:contiofparti}
  Let $\alpha\in\mathcal{P}( M^f)$ and $\{\zeta_n\}$ be a sequence of increasing measurable {partitions of $ M^f$} with $\zeta_n\nearrow\zeta$. Then
  \begin{enumerate}[label=(\roman*)]
    \item $\lim_{n\to\infty}I_{\tilde{\mu}}(\alpha|\zeta_n)(\tilde{x})=I_{\tilde{\mu}}(\alpha|\zeta)(\tilde{x})$ for ${\mu}$-a.e. $\tilde{x}\in M^f$, and
    \item\label{lem:contiofpartiitem2} $\lim_{n\to\infty}H_{\tilde{\mu}}(\alpha|\zeta_n)=H_{\tilde{\mu}}(\alpha|\zeta)$.
  \end{enumerate}
\end{lem}
\subsection{Equivalence of two definitions of unstable metric entropy}
In this subsection, we give the proof of Theorem \ref{thm:localvsfinite}, that is, we prove that the two definitions of unstable metric entropy are equivalent when $\tilde{\mu}$ is ergodic. The proof involves the relationship between two types of measurable partitions, $\eta$ and $\xi$, where the first one is a measurable partition subordinate to $W^u$-foliation constructed as in Section \ref{sec:pre}, and the latter is an increasing measurable partition subordinate to $W^u$-foliation as in Proposition \ref{prop:specialpartition}.
\begin{potA}\rm
The proof will be divided into five steps.
\begin{enumerate}[fullwidth,listparindent=2em,label=\textbf{Step \arabic*.}]
  \item\label{step:1} In this step, we show that $h^u_\mu(f,\alpha|\eta)$ is independent of $\eta$ and $\alpha$.  Firstly, let us show that for $\eta_1$ and $\eta_2\in\mathcal{P}^u( M^f)$, we have
  \[
    h_{\mu}({f},\alpha|\eta_1)=h_{\mu}({f},\alpha|\eta_2).
  \]

  By Lemma \ref{lem:info}, we have
  \begin{align}\label{eq:compute}
    H_{\tilde{\mu}}(\alpha_0^{n-1}|\eta_1)+H_{\tilde{\mu}}(\eta_2|\alpha_0^{n-1}\vee\eta_1)=& H_{\tilde{\mu}}(\alpha_0^{n-1}|\eta_2\vee\eta_1)+H_{\tilde{\mu}}(\eta_2|\eta_1), \notag\\
    H_{\tilde{\mu}}(\alpha_0^{n-1}|\eta_2)+H_{\tilde{\mu}}(\eta_1|\alpha_0^{n-1}\vee\eta_2)=& H_{\tilde{\mu}}(\alpha_0^{n-1}|\eta_1\vee\eta_2)+H_{\tilde{\mu}}(\eta_1|\eta_2).
  \end{align}
  By the construction of $\eta_1$ and $\eta_2$, we know that there are two finite partitions $\alpha_1$ and $\alpha_2$ such that $\eta_j(\tilde{x})=\alpha_j(\tilde{x})\cap \widetilde{W}^u_{\text{loc}}(\tilde{x})$, $j=1,2$, for all $\tilde{x}\in  M^f$. Let $N_1$ and $N_2$ be the cardinality of ${\alpha_1}$ and ${\alpha_2}$ respectively. For any $\tilde{x}\in M^f$, $\eta_1(\tilde{x})$ intersects at most $N_2$ elements of ${\alpha_2}$, hence intersects at most $N_2$ elements of $\eta_2$. Thus, we have
  \[
    \lim_{n\to\infty}\frac{1}{n}H_{\tilde{\mu}}(\eta_2|\alpha_0^{n-1}\vee\eta_1)\leq\lim_{n\to\infty}\frac{1}{n}H_{\tilde{\mu}}(\eta_2|\eta_1)\leq\lim_{n\to\infty}\frac{1}{n}\log N_2=0.
  \]

  Similarly, we have
  \[
    \lim_{n\to\infty}\frac{1}{n}H_{\tilde{\mu}}(\eta_1|\alpha_0^{n-1}\vee\eta_2)\leq
    \lim_{n\to\infty}\frac{1}{n}H_{\tilde{\mu}}(\eta_1|\eta_2)\leq\lim_{n\to\infty}\frac{1}{n}\log N_1=0.
  \]
  Hence we by \eqref{eq:compute}, we get
  \[
    \limsup_{n\to\infty}\frac{1}{n}H_{\tilde{\mu}}(\alpha_0^{n-1}|\eta_1)=\limsup_{n\to\infty}\frac{1}{n}H_{\tilde{\mu}}(\alpha_0^{n-1}|\eta_2).
  \]

  Then we show that for any $\beta$, $\gamma\in\mathcal{P}( M^f)$,
  \[
    \limsup_{n\to\infty}\frac{1}{n}H_{\tilde{\mu}}(\beta_0^{n-1}|\eta)=\limsup_{n\to\infty}\frac{1}{n}H_{\tilde{\mu}}(\gamma_0^{n-1}|\eta).
  \]
  Again, by Lemma \ref{lem:info}, we have
  \begin{equation}\label{eq:betagamma1}
    H_{\tilde{\mu}}(\beta_0^{n-1}|\eta)\leq H_{\tilde{\mu}}(\gamma_0^{n-1}|\eta)+H_{\tilde{\mu}}(\beta_0^{n-1}|\gamma_0^{n-1}\vee\eta),
  \end{equation}
  and similar to the proof of Lemma 2.7 (\romannumeral2) in  \cite{HuHuaWu2017}, we can show that
  \begin{equation}\label{eq:betagamma2}
    \lim_{n\to\infty}\frac{1}{n}H_{\tilde{\mu}}(\beta_0^{n-1}|\gamma_0^{n-1}\vee\eta)=0.
  \end{equation}
  By \eqref{eq:betagamma1} and \eqref{eq:betagamma2}, we have
  \[
    \limsup_{n\to\infty}\frac{1}{n}H_{\tilde{\mu}}(\beta_0^{n-1}|\eta)\leq \limsup_{n\to\infty}\frac{1}{n}H_{\tilde{\mu}}(\gamma_0^{n-1}|\eta).
  \]
  Interchanging $\beta$ with $\gamma$, we obtain
  \[
    \limsup_{n\to\infty}\frac{1}{n}H_{\tilde{\mu}}(\beta_0^{n-1}|\eta)\geq\limsup_{n\to\infty}\frac{1}{n}H_{\tilde{\mu}}(\gamma_0^{n-1}|\eta).
  \]
  \item\label{step:2} In this step, we present a construction of the increasing partition $\xi$, which will be crucial in subsequent steps. The reader can refer to Section {\uppercase\expandafter{\romannumeral9}}.2.2 in \cite{QianXieZhu2009} for more details. Given an ergodic $\tilde\mu\in\mathcal{M}({\tau})$, we can choose a set $\tilde{\Lambda}\subset M^f$, $\tilde{x}_*\in\tilde{\Lambda}$ and positive constants $\hat{\epsilon}$, $\hat{r}$ such that
\[
  B_{\tilde{\Lambda}}:=B_{\tilde{\Lambda}}(\tilde{x}_*,\hat{\epsilon}\hat{r}/2)=\{\tilde{x}\in M^f\colon \tilde{d}(\tilde{x},\tilde{x}_*)<\hat{\epsilon}\hat{r}/2\}
\]
has positive ${\tilde\mu}$ measure and the following construction of a partition $\xi_u$ satisfies Proposition \ref{prop:specialpartition}.

For each $r\in[\hat{r}/2,\hat{r}]$, put
\[
  S_{u,r}=\bigcup_{\tilde{x}\in B_{\tilde{\Lambda}}}S_u(\tilde{x},r),
\]
where $S_u(\tilde{x},r)=\{\tilde{y}\in\widetilde{W}^u_{\text{loc}}(\tilde{x})|\Pi(\tilde{y})\in B(\Pi(\tilde{x}_*),r)\}$. Then we can define a partition $\hat{\xi}_{u,\tilde{x}_*}$ of $ M^f$ such that
\[
  (\hat{\xi}_{u,\tilde{x}_*})(\tilde{y})=
  \begin{cases}
    S_u(\tilde{x},r), & {\tilde{y}\in S_u(\tilde{x},r) \text{\ for some\ }\tilde{x}\in B_{\tilde{\Lambda}},} \\
    M^f\setminus S_{u,r}, & \mbox{otherwise}.
  \end{cases}
\]

Next we can choose an appropriate $r\in[\hat{r}/2,\hat{r}]$ such that
\[
  \xi_u=\bigvee_{j=0}^\infty\tau^j\hat{\xi}_{u,\tilde{x}_*}
\]
is subordinate to $W^u$-foliation. The notation $\hat{\xi}_{u,-k}=\bigvee_{j=0}^k\tau^j\hat{\xi}_{u,\tilde{x}_*}$ will be used in the following steps.
  \item In this step, some facts concerning $\hat{\xi}_{u,-k}$ will be given, which are useful for the proof of our results.
  \begin{fact}\label{fact:convergence}
    Let $\tilde\mu\in\mathcal{M}(\tau)$ be an ergodic measure. Suppose $\eta\in\mathcal{P}^u( M^f)$ is subordinate to $W^u$-foliation, and $\hat{\xi}_{u,-k}$ is a partition described in \ref{step:2}, where $k\in\mathbb{N}\cup\{\infty\}$. Then for $\tilde\mu$-almost every ${\tilde{x}}$, there exists $N=N(\tilde{x})>0$ such that for any $j>N$, we have
  \[
    (\hat{\xi}_{u,-k-j}\vee\tau^j\eta)(\tau^j{\tilde{x}})=(\hat{\xi}_{u,-k-j})(\tau^j{\tilde{x}}).
  \]
  Hence, for any partition $\beta$ of $ M^f$ with $H_{\tilde{\mu}}(\beta|\hat{\xi}_{u,-k})<\infty$,
  \[
    I_{\tilde{\mu}}(\beta|\hat{\xi}_{u,-k-j}\vee\tau^j\eta)(\tau^j{\tilde{x}})=I_{\tilde{\mu}}(\beta|\hat{\xi}_{u,-k-j})(\tau^j{\tilde{x}}),
  \]
  which implies that
  \[
    \lim_{j\to\infty}H_{\tilde{\mu}}(\beta|\hat{\xi}_{u,-k-j}\vee\tau^j\eta)=H_{\tilde{\mu}}(\beta|\xi_u).
  \]
  Particularly, if we take $k=\infty$, then the above two equalities become
 \[
    I_{\tilde{\mu}}(\beta|\xi_u\vee\tau^j\eta)(\tau^j{\tilde{x}})=I_{\tilde{\mu}}(\beta|\xi_u)(\tau^j{\tilde{x}}),
  \]
  and
  \[
    \lim_{j\to\infty}H_{\tilde{\mu}}(\beta|\xi_u\vee\tau^j\eta)=H_{\tilde{\mu}}(\beta|\xi_u).
  \]
  \end{fact}
  \begin{pof1}\rm
  Define $\widetilde{B}(\tilde{x},\rho)$ as follows:
  \[
    \widetilde{B}^u(\tilde{x},\rho)\colon=\{\tilde{y}\in\widetilde{W}^u_{\text{loc}}(\tilde{x},f)\colon d^u_{\tilde{x}}(\Pi(\tilde{y}),\Pi(\tilde{x}))<\rho\}.
  \]

  Since $\eta$ is subordinate to $W^u$, for ${\tilde\mu}$-a.e. $\tilde{x}$, there exists $\rho=\rho(\tilde{x})>0$ such that $\widetilde{B}^u(\tilde{x},\rho)\subset\eta(\tilde{x})$. Since ${\tilde\mu}$ is ergodic, for ${\tilde\mu}$-a.e. $\tilde{x}\in M^f$, there are infinitely many $n>0$ such that $\tau^n\tilde{x}\in S_{u,r}$. Take $N=N(\tilde{x})$ large enough such that
  \[
    \tau^N\tilde{x}\in S_{u,r}
  \]
  and
  \[
   \tau^{-N}(\hat{\xi}_{u,\tilde{x}_*}(\tau^N\tilde{x}))\subset \widetilde{B}^u(\tilde{x},\rho)\subset\eta(\tilde{x}).
  \]
  Then we have
  \[
    {\tau^{-j}}\Big(\tau^{j-N}(\hat{\xi}_{u,{\tilde{x}_*}}({\tau^N\tilde{x}}))\Big)\subset\eta(\tilde{x})
  \]
  for any $j\geq N$. Since
  \[
    \hat{\xi}_{u,-k-j}=\bigvee_{l=0}^{k+j}\tau^l(\hat{\xi}_{u,{\tilde{x}_*}})\geq \tau^{j-N}(\hat{\xi}_{u,{\tilde{x}_*}}),
  \]
  so we have
  \[
    {\tau^{-j}}\Big((\hat{\xi}_{u,-k-j})(\tau^j(\tilde{x}))\Big)\subset\eta(\tilde{x}).
  \]
  Thus we have
  \[
    (\hat{\xi}_{u,-k-j})(\tau^j(\tilde{x}))\subset(\tau^j\eta)(\tau^j(\tilde{x})),
  \]
  which implies that
  \[
    \left((\hat{\xi}_{u,-k-j})\vee \tau^j\eta\right)(\tau^j(\tilde{x}))=(\hat{\xi}_{u,-k-j})(\tau^j(\tilde{x})).
  \]
 This proves the first statement in Fact \ref{fact:convergence}.

  Then following the line of the proof for Lemma 2.11 in  \cite{HuHuaWu2017}, we can prove the remaining results in Fact \ref{fact:convergence}, where Fatou's Lemma and Lemma \ref{lem:contiofparti} are needed.
  \end{pof1}

  The proof of the following fact is analogous to that in  \cite{HuHuaWu2017}, the reader can refer to the proof of Lemma 2.10 in  \cite{HuHuaWu2017} for more details.
  \begin{fact}\label{fact:approach}
    Suppose that $\tilde{\mu}\in\mathcal{M}(\tau)$ is an ergodic measure and $\alpha\in\mathcal{P}( M^f)$ is finite. For any $\epsilon>0$, there exists $K>0$ such that for any $k\geq K$,
  \[
    \limsup_{n\to\infty}H_{\tilde{\mu}}(\alpha|\alpha^n_1\vee(\hat{\xi}_{u,-k})_1^n)\leq\epsilon.
  \]
  \end{fact}

  The following fact comes from  \cite{QianXieZhu2009}.
  \begin{fact}[Proposition \uppercase\expandafter{\romannumeral9}.3.1 in  \cite{QianXieZhu2009}]\label{fact:3}
    When $\tilde{\mu}\in\mathcal{M}(\tau)$ is an ergodic measure, we have
    \[
      \tilde{h}^u_\mu(f)=H_{\tilde{\mu}}(\xi_u|\tau\xi_u).
    \]
  \end{fact}
  \item In this step, we prove that $h_\mu(f,\alpha|\eta)\leq\tilde{h}^u_\mu(f)$. By Lemma \ref{lem:info2} \ref{lem:info2_item1}, with $\gamma=\eta$ and $\beta=\hat{\xi}_{u,-k}$, we have for any $\eta\in\mathcal{P}^u( M^f)$, $n>0$,
  \begin{equation}\label{eq:step4}
    \frac{1}{n}H_{\tilde{\mu}}((\hat{\xi}_{u,-k})_0^{n-1}|\eta)=\frac{1}{n}H_{\tilde{\mu}}(\hat{\xi}_{u,-k}|\eta)+\frac{1}{n}\sum_{j=0}^{n-1}H_{\tilde{\mu}}(\hat{\xi}_{u,-k}|\tau\hat{\xi}_{u,-k-j+1}\vee\tau^j\eta).
  \end{equation}
  By Fact \ref{fact:convergence}, the second term of the right side of \eqref{eq:step4} converges to $H_{\tilde{\mu}}(\hat{\xi}_{u,-k}|\tau\xi_u)$ as $j\to\infty$. It is clear that each elements of $\eta$ intersects at most {$2^{k+1}$} elements of $\hat{\xi}_{u,-k}$. So we have
  \[
    H_{\tilde{\mu}}(\hat{\xi}_{u,-k}|\eta)\leq\log{2^{k+1}},
  \]
  which implies that
  \[
    \lim_{n\to\infty}\frac{1}{n}H_{\tilde{\mu}}(\hat{\xi}_{u,-k}|\eta)=0.
  \]
  Thus we get
  \begin{equation}\label{eq:finitevsmeas}
    \lim_{n\to\infty}\frac{1}{n}H_{\tilde{\mu}}((\hat{\xi}_{u,-k})_0^{n-1}|\eta)=H_{\tilde{\mu}}(\hat{\xi}_{u,-k}|\tau\xi_u)\leq H_{\tilde{\mu}}(\xi_u|\tau\xi_u).
  \end{equation}

  By Lemma \ref{lem:info2} \ref{lem:info2_item2} with $\gamma=(\hat{\xi}_{u,-k})_0^{n-1}$ and the fact that
  \[
    \tau^j(\hat{\xi}_{u,-k})_0^{n-1}=(\hat{\xi}_{u,-k-j})_0^{n-j-1},
  \]
  we know that
  \begin{align*}
    H_{\tilde{\mu}}(\alpha_0^{n-1}|(\hat{\xi}_{u,-k})_0^{n-1}) &=H_{\tilde{\mu}}(\alpha|\hat{\xi}_{u,-n-k+1})+\sum_{j=0}^{n-2}H_{\tilde{\mu}}(\alpha|\alpha_1^{n-1-j}\vee(\hat{\xi}_{u,-k-j})_0^{n-1-j}) \\
    &=H_{\tilde{\mu}}(\alpha|\hat{\xi}_{u,-n-k+1})+\sum_{j=1}^{n-1}H_{\tilde{\mu}}(\alpha|\alpha_1^{j}\vee{\hat{\xi}^j}_{u,-k-n+1+j}) \\
    &\leq H_{\tilde{\mu}}(\alpha)+\sum_{j=1}^{n-1}H_{\tilde{\mu}}(\alpha|\alpha_1^j\vee(\hat{\xi}_{u,-k})_1^j).
  \end{align*}

  For any $\epsilon>0$, take $k>0$ as in Fact \ref{fact:approach}, thus we have
  \[
    \limsup_{n\to\infty}H_{\tilde{\mu}}(\alpha|\alpha_1^{n-1}\vee(\hat{\xi}_{u,-k})_1^{n-1})\leq \epsilon.
  \]
  Then we get
   \begin{equation}\label{eq:finitevsmeas2}
    \limsup_{n\to\infty}\frac{1}{n}H_{\tilde{\mu}}(\alpha_0^{n-1}|(\hat{\xi}_{u,-k})_0^{n-1})\leq\epsilon.
  \end{equation}

  By Lemma \ref{lem:info}, we have
  \begin{equation}\label{eq:finitevsmeas3}
    H_{\tilde{\mu}}(\alpha_0^{n-1}|\eta)\leq H_{\tilde{\mu}}((\hat{\xi}_{u,-k})_0^{n-1}|\eta)+H_{\tilde{\mu}}(\alpha_0^{n-1}|(\hat{\xi}_{u,-k})^{n-1}_0).
  \end{equation}

  Thus, by \eqref{eq:finitevsmeas2}, \eqref{eq:finitevsmeas3}, then by \eqref{eq:finitevsmeas} and Fact \ref{fact:3} we have
  \begin{align*}
    h_\mu({f},\alpha|\eta) &=\limsup_{n\to\infty}\frac{1}{n}H_{\tilde{\mu}}(\alpha_0^{n-1}|\eta) \\
    &\leq \lim_{n\to\infty}\frac{1}{n}H_{\tilde{\mu}}((\hat{\xi}_{u,-k})_0^{n-1}|\eta)+\epsilon \\
    &{\leq H_{\tilde{\mu}}(\xi_u|\tau\xi_u)+\epsilon} \\
    &=\tilde{h}^u_\mu(f)+\epsilon.
  \end{align*}
  Let $\epsilon\to0$, we have $h_\mu({f},\alpha|\eta)\leq \tilde{h}^u_\mu(f).$
  \item In this step, we complete the proof of Theorem \ref{thm:localvsfinite}. By a similar treatment in the proof of Proposition 2.13 in  \cite{HuHuaWu2017}, we can construct an increasing measurable partition $\tilde{\xi}$ satisfying Proposition \ref{prop:specialpartition} with diameter bounded above. And we know that $h_\mu({f},\tilde{\xi})=h_\mu({f},\xi_u)$.
  So we only need to prove $h_\mu(f,\alpha|\eta)\leq\tilde{h}^u_\mu(f)$ for $\tilde{\xi}$.

  We can choose a sequence of partitions $\alpha_n\in\mathcal{P}( M^f)$ such that
  \[
    \mathcal{B}(\alpha_n)\nearrow\mathcal{B}(\tau^{-1}\tilde{\xi})\text{ as }n\to\infty,
  \]
  which implies
  \[
    \lim_{n\to\infty}H_{\tilde{\mu}}(\alpha_n|\tilde{\xi})=H_{\tilde{\mu}}(\tau^{-1}\tilde{\xi}|\tilde{\xi}).
  \]
  Thus, we have
  \[
    \sup_{\alpha\in\mathcal{P}( M^f),\alpha<\tau^{-1}\tilde{\xi}}H_{\tilde{\mu}}(\alpha|\tilde{\xi})=H_{\tilde{\mu}}(\tau^{-1}\tilde{\xi}|\tilde{\xi}).
  \]

  For any $\alpha\in\mathcal{P}( M^f)$ with $\alpha<\tau^{-1}\tilde{\xi}$, we have that for any $j>0$, $\tau^j\alpha^{j-1}_0<\tau^j(\tau^{-1}\tilde{\xi})^{j-1}_0=\tilde{\xi}$. Then by Lemma \ref{lem:info2} \ref{lem:info2_item1}, we have
  \begin{align*}
    H_{\tilde{\mu}}(\alpha_0^{n-1}|\eta) &= H_{\tilde{\mu}}(\alpha|\eta)+\sum_{j=1}^{n-1}H_{\tilde{\mu}}(\alpha|\tau^j(\alpha_0^{j-1}\vee\eta)) \\
    &\geq H_{\tilde{\mu}}(\alpha|\eta)+\sum_{j=1}^{n-1}H_{\tilde{\mu}}(\alpha|\tilde{\xi}\vee\tau^j\eta).
  \end{align*}
  Then by Fact \ref{fact:convergence} we have
  \[
    \lim_{j\to\infty}H_{\tilde{\mu}}(\alpha|\tilde{\xi}\vee\tau^j\eta)=H_{\tilde{\mu}}(\alpha|\tilde{\xi}),
  \]
  which implies that
  \[
    \limsup_{n\to\infty}\frac{1}{n}H_{\tilde{\mu}}(\alpha_0^{n-1}|\eta)\geq\liminf_{n\to\infty}\frac{1}{n}H_{\tilde{\mu}}(\alpha_0^{n-1}|\eta)\geq H_{\tilde{\mu}}(\alpha|\tilde{\xi}).
  \]
  So we have
  \begin{align*}
    \sup_{\alpha\in\mathcal{P}( M^f)}h_\mu({f},\alpha|\eta) &\geq\sup_{\alpha\in\mathcal{P}( M^f),\alpha<\tau^{-1}\tilde{\xi}}h_\mu({f},\alpha|\eta) \\
    &= \sup_{\alpha\in\mathcal{P}( M^f),\alpha<\tau^{-1}\tilde{\xi}}\limsup_{n\to\infty}\frac{1}{n}H_{\tilde{\mu}}(\alpha_0^{n-1}|\eta) \\
    &\geq \sup_{\alpha\in\mathcal{P}( M^f),\alpha<\tau^{-1}\tilde{\xi}}\liminf_{n\to\infty}\frac{1}{n}H_{\tilde{\mu}}(\alpha_0^{n-1}|\eta) \\
    &\geq \sup_{\alpha\in\mathcal{P}( M^f),\alpha<\tau^{-1}\tilde{\xi}}H_{\tilde{\mu}}(\alpha|\tilde{\xi}) \\
    &=H_{\tilde{\mu}}(\tau^{-1}\tilde{\xi}|\tilde{\xi}).
  \end{align*}
  By the statement in \ref{step:1}, $h_\mu({f},\alpha|\eta)$ is independent of $\alpha$, meaning
  \[
    h_\mu({f},\alpha|\eta)=\sup_{\beta\in\mathcal{P}( M^f)}h_\mu({f},\beta|\eta)
  \]
    for any $\alpha\in\mathcal{P}( M^f)$. This finishes the proof of Theorem \ref{thm:localvsfinite}.
\end{enumerate}
\end{potA}

A corollary can be obtained directly as follows.
\begin{cor}
  Suppose that $\tilde\mu\in\mathcal{M}(\tau)$ is ergodic, then for any $\alpha\in\mathcal{P}( M^f)$ and $\eta\in\mathcal{P}^u( M^f)$, we have
  \[
    {h}_\mu^u({f})=h_\mu({f},\alpha|\eta)=\lim_{n\to\infty}\frac{1}{n}H_{\tilde{\mu}}(\alpha^{n-1}_0|\eta).
  \]
\end{cor}

\subsection{Shannon-McMillan-Breiman Theorem for unstable metric entropy}
In this subsection, we give the proof of Theorem \ref{thm:SMB}. We always assume that $\tilde{\mu}\in\mathcal{M}(\tau)$ is ergodic. Firstly, we need some lemmas.
\begin{lem}\label{lem:SMBlow}
  Let $\alpha\in\mathcal{P}( M^f)$, $\eta\in\mathcal{P}^u( M^f)$. Then for any $\xi\in\mathcal{Q}^u( M^f)$, we have
  \[
     {h_\mu({f},\alpha|\eta)}\leq\liminf_{n\to\infty}\frac{1}{n}I_{\tilde{\mu}}(\alpha_0^{n-1}|\xi)(\tilde{x})\quad\text{ for }{\tilde\mu}\text{-a.e. }\tilde{x}.
  \]
\end{lem}

\begin{pf}\rm
For any $\epsilon>0$, there exists $k>0$ such that $\mathrm{diam}(\alpha^k_0\vee\xi)\leq\epsilon$. Then for $n>0$, we have
  \[
    (\alpha_0^{k+n-1}\vee\xi)(\tilde{x})=\bigvee_{j=0}^{n-1}( {\tau^{-j}}\alpha^k_0\vee\xi)(\tilde{x})\subset \widetilde{V}^u({f},\tilde{x},n,\epsilon).
  \]
  By Theorem \ref{thm:localvsfinite} and Remark \ref{rmk:low=up} we know that for ${\tilde\mu}$-a.e. $\tilde{x}$,
  \begin{align*}
    h_\mu({f},\alpha|\eta)&=h_\mu({f},\tilde{x},\xi)\\
   & = \lim_{n\to\infty}-\frac{1}{n}\log{\mu}_{\tilde{x}}^{\xi}\widetilde{V}^u({f},\tilde{x},n,\epsilon) \\
    &\leq \liminf_{n\to\infty}-\frac{1}{n}\log{\mu}_{\tilde{x}}^{\xi}((\alpha^{k+n-1}_0)(\tilde{x})) \\
    &= \liminf_{n\to\infty}-\frac{1}{n}\log{\mu}_{\tilde{x}}^{\xi}((\alpha^{n-1}_0)(\tilde{x})) \\
    &= \liminf_{n\to\infty}\frac{1}{n}I_{\tilde{\mu}}(\alpha^{n-1}_0|\xi)(\tilde{x}).
  \end{align*}
  \end{pf}

The following lemmas are counterparts of those in  \cite{HuHuaWu2017}, which are completely parallel to the treatment in  \cite{HuHuaWu2017}, so we omit the proofs.
\begin{lem}[Lemma 3.4 in  \cite{HuHuaWu2017}]\label{lem:SMBlow2}
  Let $\eta\in\mathcal{P}^u( M^f)$ and $\xi\in\mathcal{Q}^u( M^f)$. Then for ${\tilde\mu}$-a.e. $\tilde{x}$, we have
  \[
    \liminf_{n\to\infty}\frac{1}{n}I_{\tilde{\mu}}(\alpha_0^{n-1}|\xi)(\tilde{x})=\liminf_{n\to\infty}\frac{1}{n}I_{\tilde{\mu}}(\alpha_0^{n-1}|\eta)(\tilde{x}),
  \]
  \[
    \limsup_{n\to\infty}\frac{1}{n}I_{\tilde{\mu}}(\alpha_0^{n-1}|\xi)(\tilde{x})=\limsup_{n\to\infty}\frac{1}{n}I_{\tilde{\mu}}(\alpha_0^{n-1}|\eta)(\tilde{x}).
  \]
\end{lem}
\begin{lem}[Lemma 3.7 in  \cite{HuHuaWu2017}]\label{lem:SMBup}
  For any $\eta\in\mathcal{P}^u( M^f)$ and $\xi\in\mathcal{Q}^u( M^f)$, we have
  \[
    \lim_{n\to\infty}\frac{1}{n}I_{\tilde{\mu}}(\tau^{-n}\xi|\eta)(\tilde{x})=\lim_{n\to\infty}\frac{1}{n}I_{\tilde{\mu}}(\tau^{-n}\xi|\xi)(\tilde{x})=h_\mu({f},\tilde{x},\xi).
  \]
\end{lem}
\begin{lem}[Lemma 3.8 in  \cite{HuHuaWu2017}])\label{lem:SMBup2}
  Let $\alpha\in\mathcal{P}( M^f)$, $\eta\in\mathcal{P}^u( M^f)$. Then for ${\tilde\mu}$-a.e. $\tilde{x}$, we have
  \[
    \lim_{n\to\infty}\frac{1}{n}I_{\tilde{\mu}}(\alpha^{n-1}_0|\xi^{n-1}_0\vee\eta)(\tilde{x})=0.
  \]
\end{lem}

Now, we begin to prove Theorem \ref{thm:SMB}.

\begin{potB}\rm
  By Lemma \ref{lem:SMBlow} and Lemma \ref{lem:SMBlow2} we can get directly
  \begin{equation}\label{eq:SMBlow}
    h_\mu({f},\alpha|\eta)\leq\liminf_{n\to\infty}\frac{1}{n}I_{\tilde{\mu}}(\alpha_0^{n-1}|\eta)(\tilde{x}).
  \end{equation}
  By Lemma \ref{lem:info}, we have
  \begin{align*}
    I_{\tilde{\mu}}(\alpha^{n-1}_0|\eta)(\tilde{x}) \leq& I_{\tilde{\mu}}(\alpha^{n-1}_0\vee\xi^{n-1}_0|\eta)(\tilde{x}) \\
   =& I_{\tilde{\mu}}(\xi^{n-1}_0|\eta)(\tilde{x})+I_{\tilde{\mu}}(\alpha^{n-1}_0|\xi^{n-1}_0\vee\eta)(\tilde{x}).
  \end{align*}
  Then by Lemma \ref{lem:SMBup2}, Lemma \ref{lem:SMBup}, and Theorem \ref{thm:localvsfinite}, we have
  \begin{align}\label{eq:SMBup}
    \limsup_{n\to\infty}\frac{1}{n}I_{\tilde{\mu}}(\alpha^{n-1}_0|\eta)(\tilde{x}) &\leq\limsup_{n\to\infty}\frac{1}{n}I_{\tilde{\mu}}(\xi^{n-1}_0|\eta)(\tilde{x}) \notag \\
    &= h_\mu^u({f})=h_\mu({f},\alpha|\eta).
  \end{align}
  Combining \eqref{eq:SMBlow} and \eqref{eq:SMBup}, we complete the proof of Theorem \ref{thm:SMB}.
\end{potB}
\section{Unstable topological entropy and unstable pressure for endomorphisms}\label{sec:pressure}
In this section, we give the definition of unstable topological entropy and unstable pressure for a potential function $\varphi\in C(M)$ for endomorphisms.

Similar to the classical pressure, there are several ways to define unstable pressure. Here we use $W^u$-separated sets. Fix $\delta>0$, for $\tilde{x}\in M^f$, Let $\overline{W^u(\tilde{x},\delta)}$ be the $\delta$-neighborhood of $x_0$ in $W^u(\tilde{x},f)$. A subset $E$ of $\overline{W^u(\tilde{x},\delta)}$ is called an $(n,\epsilon)$ {$W^u$-separated set} if for any $y_1,y_2\in E$, we have $d^u_{n}(y_1,y_2)>\epsilon$, where $ d^u_{n}(y_1,y_2)$ is defined by
{\[
  d^u_{n}(y_1,y_2):=\max_{0\leq j\leq n-1}\{ d^u_{\tau^j\tilde{x}}(f^j(y_1),f^j(y_2))\}.
\]}
Recall that $d^u_{\tilde{x}}$ is the metric on $W^u(\tilde{x},f)$ induced by the Riemannian structure on $M$.

Now we can define $P^u(f,\varphi,\tilde{x},\delta,n,\epsilon)$ as follows,
\begin{align*}
 P^u(f,\varphi,\tilde{x},\delta,n,\epsilon)=\sup & \Big\{\sum_{y\in E}\exp((S_n\varphi)(y)):\\ & \text{ } E\text{ is an }(n,\epsilon)\ W^u\text{-separated set of }\overline{W^u(\tilde{x},\delta)}\Big\},
\end{align*}
{where $(S_n\varphi)(y)=\sum_{j=0}^{n-1}\varphi(f^j(y))$.} Then $P^u(f,\varphi,\tilde{x},\delta)$ is defined as
\[
  P^u(f,\varphi,\tilde{x},\delta)=\lim_{\epsilon\to 0}\limsup_{n\to\infty}\frac{1}{n}\log P^u(f,\varphi,\tilde{x},\delta,n,\epsilon).
\]
Next, we define
\[
  P^u(f,\varphi,\delta)=\sup_{\tilde{x}\in M^f}P^u(f,\varphi,\tilde{x},\delta)
\]

Let $\tilde{\varphi}(\tilde{x})=\varphi(\Pi(\tilde{x}))$. It is easy to check that
\[
  \int_{M^f}\tilde{\varphi}\mathrm{d}\tilde{\mu}=\int_M\varphi\mathrm{d}\mu.
\]
Denote
\[
  \textcolor[rgb]{0.00,0.00,0.00}{\overline{\widetilde{W}^u(\tilde{x},\delta)}=\{\tilde{y}\in M^f\colon\Pi(\tilde{y})\in \overline{W^u(\tilde{x},\delta)}\text{ and satisfies \eqref{eq:locunstable2}}\}.}
\]
A subset $\widetilde{E}$ of $\overline{\widetilde{W}^u(\tilde{x},\delta)}$ is called an $(n,\epsilon)$ $W^u$-separated set if for any $\tilde{y}_1$, $\tilde{y}_2\in\widetilde{E}$, we have
\[
  \tilde{d}^u_{n}(\tilde{y}_1,\tilde{y}_2)>\epsilon.
\]
Then we can define
\begin{align*}
 \widetilde{P}(\tau,{\tilde{\varphi}},\tilde{x},\delta,n,\epsilon)\colon=\sup & \Big\{\sum_{\tilde{y}\in \widetilde{E}}\exp((\tilde{S}_n\tilde{\varphi})(\tilde{y})):\\ & \text{ } \widetilde{E}\text{ is an }(n,\epsilon)\ W^u\text{-separated set of }\overline{\widetilde{W}^u(\tilde{x},\delta)}\Big\},
\end{align*}
{where $(\tilde{S}_n\varphi)(\tilde{y})=\sum_{j=0}^{n-1}\tilde\varphi(\tau^j(\tilde{y}))$.}

It is clear that for an $(n,\delta)$ $W^u$-separated set $\widetilde{E}$ of $\overline{\widetilde{W}^u(\tilde{x},\delta)}$, there is an $(n,\delta)$ $W^u$-separated set $E$ with the same cardinality as $\widetilde{E}$, and vice versa. Then noticing that $\varphi(f^j(\Pi(\tilde{x})))=\tilde{\varphi}(\tau^j(\tilde{x}))$ we have $\widetilde{P}(\tau,{\tilde{\varphi}},\tilde{x},\delta,n,\epsilon)=P^u(f,\varphi,\tilde{x},\delta,n,\epsilon)$. Then $\widetilde{P}(\tau,\tilde{\varphi},\tilde{x},\delta)$ and $\widetilde{P}(\tau,\tilde{\varphi},\delta)$ can be formulated similarly.

Finally, we can give the definition of unstable pressure for $f$.
\begin{mdef}\rm
  The \textit{unstable pressure} for $f$ is defined as
  \[
    P^u(f,\varphi)=\lim_{\delta\to 0}P^u(f,\varphi,\delta)=\lim_{\delta\to 0}\widetilde{P}(\tau,\tilde{\varphi},\delta).
  \]
\end{mdef}

We can prove that $P^u(f,\varphi)$ \textit{is independent of $\delta>0$}. Indeed, notice that for given $\delta_1<\delta$ and $\tilde{x}\in M^f$, there exists a positive number $N=N(\delta,\delta_1)$ depending on the Riemannian structure on $\overline{\widetilde{W}^u(\tilde{x},\delta)}$ such that
\[
  \overline{\widetilde{W}^u(\tilde{x},\delta)}\subset\bigcup_{j=1}^{N}\overline{\widetilde{W}^u(\tilde{y}_j,\delta_1)}
\]
for some $\tilde{y}_j\in\overline{\widetilde{W}^u(\tilde{x},\delta)}$, $j=1,2,\cdots,N$. Then following the calculation in the proof of Lemma 4.1 in  \cite{HuHuaWu2017}, we can prove that $P^u({f},\varphi,\delta)\leq P^u({f},\varphi)$, and it is clear that $P^u({f},\varphi,\delta)\geq P^u({f},\varphi)$, which means $P^u(f,\varphi)$ {is independent of $\delta$}.
\begin{mdef}\rm
  The \textit{unstable topological entropy} of $f$ is defined as
  \[
    h^u_{\text{{\rm top}}}(f)=P^u(f,0).
  \]
\end{mdef}

The following proposition can be obtained directly from the definitions.
\begin{prop}
  For any $\varphi$, $\psi\in C(M)$ and constant $c\in \mathbb{R}$, the following properties hold.
  \begin{enumerate}[label=(\roman*)]
    \item If $\varphi\leq \psi$, then $P^u(f ,\varphi)\leq P^u(f ,\psi)$;
    \item $P^u(f ,\varphi+c)=P^u(f ,\varphi)+c$;
    \item {$h_{\text{{\rm top}}}^u(f)+\inf \varphi \leq P^u(f, \varphi) \leq h_{\text{{\rm top}}}^u(f)+\sup\varphi$;}
    \item {if $P^u(f ,\cdot)<\infty$, $|P^u(f, \varphi)-P^u(f, \psi)|\leq \|\varphi-\psi\|$;}
    \item {if $P^u(f ,\cdot)<\infty$}, then the map $P^u(f ,\cdot)\colon C(M)\to\mathbb{R}\cup\{\infty\}$ is convex;
    \item $P^u(f ,\varphi+\psi\circ f-\psi)=P^u(f ,\varphi)$;
    \item $P^u(f ,\varphi+\psi)\leq P^u(f ,\varphi)+P^u(f ,\psi)$.
  \end{enumerate}
\end{prop}
\section{Variational principle}\label{sec:vp}
\subsection{Variational principle for unstable pressure}
In this subsection, we prove our main result of this paper, i.e. Theorem \ref{thm:vp}, whose proof consists of two parts.
\begin{potC}\rm
Let $\varphi\in C(M)$.
\begin{enumerate}[fullwidth,listparindent=2em,label=\textbf{Part \Roman*.}]
\item In this part, we prove that For $\mu\in\mathcal{M}(f)$,
  \[
    h^u_\mu(f)+\int_{M}\varphi \mathrm{d}{\mu}\leq P^u(f,\varphi).
  \]

Firstly, we give a useful lemma from  \cite{HuWuZhu2017}.
\begin{lem}\label{lem:wellknown}
  Suppose $0\leq p_1$, $\cdots$, $p_m\leq1$, {$s=p_1+\cdots+p_m$} and $a_1$, $\cdots$, $a_m\in\mathbb{R}$. Then
  \[
    \sum_{i=1}^{m}p_i(a_i-\log p_i)\leq s\left(\log\sum_{i=1}^{m}e^{a_i}-\log s\right).
  \]
\end{lem}

The following two lemmas are also important, whose proofs are analogous to those in  \cite{HuHuaWu2017}.
\begin{lem}[Proposition 2.14 in \cite{HuHuaWu2017}]\label{lem:affine}
  For any $\alpha\in\mathcal{P}( M^f)$ and $\eta\in\mathcal{P}^u( M^f)$, the map ${{\tilde{\mu}}}\mapsto H_{{{\tilde{\mu}}}}(\alpha|\eta)$ from $\mathcal{M}({\tau})$ to $\mathbb{R}^+\cup \{0\}$ is concave. Moreover, the map ${\tilde{\mu}}\mapsto h_{{\mu}}^u({f})$ from $\mathcal{M}({\tau})$ to $\mathbb{R}^+\cup \{0\}$ is affine.
\end{lem}

\begin{lem}[Proposition 2.15 in \cite{HuHuaWu2017}]\label{lem:semiconti}
  Let ${\tilde{\mu}}\in\mathcal{M}({\tau})$ and {$\eta\in \mathcal{P}^u( M^f)$.} Assume that there exists a sequence of partitions $\{\beta_n\}_{n=1}^{\infty}\subset\mathcal{P}( M^f)$ such that $\beta_1<\beta_2<\cdots<\beta_n<\cdots$ and {$\mathcal{B}(\beta_n)\nearrow\mathcal{B}(\eta)$}, {and moreover, ${{{\mu}}}(\partial(\Pi(\beta_n)))=0$, for $n=1,2,\cdots$. Let $\alpha\in\mathcal{P}( M^f)$ satisfy ${{{\mu}}}(\partial(\Pi(\alpha)))=0$}. Then the function ${{\tilde{\mu}}}'\mapsto H_{{{\tilde{\mu}}}'}(\alpha|\eta)$ is upper semi-continuous at ${{\tilde{\mu}}}$, i.e.,
  \[
    \limsup_{{{\tilde{\mu}}}'\to{{\tilde{\mu}}}}H_{{{\tilde{\mu}}}'}(\alpha|\eta)\leq H_{{{\tilde{\mu}}}}(\alpha|\eta).
  \]
  Moreover, the function ${\tilde{\mu}}'\mapsto h_{{{\mu}}'}^u({f})$ is upper semi-continuous at ${\tilde{\mu}}$, i.e.,
  \[
    \limsup_{{{\tilde{\mu}}}'\to{\tilde{\mu}}}h_{{{\mu}}'}^u({f})\leq h_{{{\mu}}}^u({f}).
  \]
\end{lem}

By the definition of unstable pressure and $\tilde{\varphi}$, we only need to prove that
\begin{equation}\label{eq:vpleq}
  h^u_\mu(f)+\int_{M^f}\tilde{\varphi}\mathrm{d}\tilde{\mu}\leq P^u(f,\varphi).
\end{equation}

Let ${\tilde{\mu}}=\int_{\mathcal{M}^e({\tau})}{\tilde{\nu}} \mathrm{d}m({\tilde{\nu}})$ be the unique ergodic decomposition
where ${\mathcal{M}^e({\tau})}$ is the set of ergodic measures in ${\mathcal{M}({\tau})}$ and $m$ is a Borel probability measure such that $m({\mathcal{M}^e({\tau})})=1$.
Since ${\tilde{\mu}} \mapsto h_{{\mu}}^u({f})$ is affine and upper semi-continuous by Lemma \ref{lem:affine} and \ref{lem:semiconti}, so is ${\tilde{\mu}} \mapsto h_{{\mu}}^u({f})+\int_{ M^f}\tilde{\varphi} \mathrm{d}{{\tilde{\mu}}}$ and hence
\begin{equation*}\label{e:ergodicdecom}
h_{{\mu}}^u({f})+\int_{ M^f}\tilde{\varphi} \mathrm{d}{{\tilde{\mu}}}=\int_{\mathcal{M}^e({\tau})}\Big(h_\nu^u({f})+\int_{ M^f}\tilde{\varphi} \mathrm{d}{{\tilde{\nu}}}\Big) \mathrm{d}m({\tilde{\nu}})
\end{equation*}
So we only need to prove \eqref{eq:vpleq} for
ergodic measures.

 We assume ${\tilde{\mu}}$ is ergodic. Let $\xi\in \mathcal{Q}^u( M^f)$, Then we can pick $\tilde{x}\in M^f$ satisfying
  \begin{enumerate}[label=(\roman*)]
    \item ${{\tilde{\mu}}}_{\tilde{x}}^\xi(\xi{(\tilde{x})})=1$;
    \item\label{prp:item2} there exists ${\widetilde{B}}\subset\xi{(\tilde{x})}$ such that
    \begin{enumerate}[label=(\alph*)]
      \item ${{\tilde{\mu}}}_{\tilde{x}}^\xi({\widetilde{B}})=1$;
      \item $h_{{\mu}}({f},\xi)=h_{{\mu}}({f},{\tilde{y}},\xi)= \lim_{n\to\infty}-\frac{1}{n}\log{{\tilde{\mu}}}^\xi_{\tilde{y}}(\widetilde{V}^u({f},{\tilde{y}},n,\epsilon))$ for any $\tilde{y}\in {\widetilde{B}}$ and $\epsilon>0$, according to Remark \ref{rmk:low=up};
      \item $\lim_{n\to\infty}\frac{1}{n}(\widetilde{S}_n\tilde{\varphi})({\tilde{y}})=\int_{ M^f}\tilde{\varphi} \mathrm{d}{{\tilde{\mu}}}$ for any $\tilde{y}\in {\widetilde{B}}$, which can be obtained using the Birkhoff ergodic theorem on $( M^f,\tau)$.
    \end{enumerate}
       \end{enumerate}

  Fix $\rho>0$. By property \ref{prp:item2} we know that for any $\tilde{y}\in {\widetilde{B}}$, there exists $N(\tilde{y})=N(\tilde{y},\epsilon)>0$ such that if $n\geq N(\tilde{y})$ then we have
  \[
    {{\tilde{\mu}}}^\xi_{\tilde{y}}(\widetilde{V}^u({f},{\tilde{y}},n,\epsilon))\leq e^{-n(h_{{\mu}}({f},\xi)-\rho)}
  \]
  and
  \begin{equation}\label{eq:est2}
    \frac{1}{n}(\widetilde{S}_n\tilde{\varphi})({\tilde{y}})\geq\int_{ M^f}\tilde{\varphi} \mathrm{d}{{\tilde{\mu}}}-\rho.
  \end{equation}

  Denote ${\widetilde{B}}_n=\{\tilde{y}\in {\widetilde{B}}\colon N(\tilde{y})\leq n\}$. Then ${\widetilde{B}}=\bigcup_{n=1}^\infty {\widetilde{B}}_n$. So we can choose $n>0$ such that ${{\tilde{\mu}}}^\xi_{\tilde{x}}({\widetilde{B}}_n)>{{\tilde{\mu}}}^\xi_{\tilde{x}}({\widetilde{B}})-\rho=1-\rho$. If $\tilde{y}\in {\widetilde{B}}_n\subset\xi{(\tilde{x})}$, then ${{\tilde{\mu}}}^\xi_{\tilde{y}}={{\tilde{\mu}}}^\xi_{\tilde{x}}$. So for any $\tilde{y}\in {\widetilde{B}}_n$ we have
  \begin{equation}\label{eq:est3}
    {{\tilde{\mu}}}^\xi_{\tilde{x}}(\widetilde{V}^u({f},{\tilde{y}},n,\epsilon))\leq e^{-n(h_{{\mu}}({f},\xi)-\rho)}.
  \end{equation}

  Now we can choose $\delta>0$ such that $\widetilde{W}^u(\tilde{x},\delta)\supset\xi{(\tilde{x})}$. Let ${\widetilde{F}}$ be an {$(n,\epsilon/2)$ $W^u$-spanning set} of $\overline{\widetilde{W}^u(\tilde{x},\delta)}\cap {\widetilde{B}}_n$ (i.e. for any $\tilde{z}\in\overline{\widetilde{W}^u(\tilde{x},\delta)}\cap {\widetilde{B}}_n$, there is $\tilde{y}\in {\widetilde{F}}$ such that $\tilde{d}^u_n(\tilde{y},\tilde{z})<\epsilon/2$.) satisfying
  \[
    \overline{\widetilde{W}^u(\tilde{x},\delta)}\cap {\widetilde{B}}_n\subset\bigcup_{\tilde{z}\in {\widetilde{F}}}\widetilde{V}^u({f},\tilde{z},n,\epsilon/2),
  \]
  and $\widetilde{V}^u({f},{\tilde{z}},n,\epsilon/2)\cap {\widetilde{B}}_n\neq\emptyset$ for any $\tilde{z}\in {\widetilde{F}}$. Then choose an arbitrary point in $\widetilde{V}^u({f},{\tilde{z}},n,\epsilon/2)\cap {\widetilde{B}}_n$, which is denoted by $\tilde{y}(\tilde{z})$. Then we have
  \begin{align}\label{eq:est3.5}
    1-\rho &< {{\tilde{\mu}}}^\xi_{\tilde{x}}(\overline{\widetilde{W}^u(\tilde{x},\delta)}\cap {\widetilde{B}}_n)\notag\\
    &\leq {{\tilde{\mu}}}^\xi_{\tilde{x}}(\bigcup_{\tilde{z}\in {\widetilde{F}}}\widetilde{V}^u({f},{\tilde{z}},n,\epsilon/2)) \notag\\
    &\leq \sum_{\tilde{z}\in {\widetilde{F}}}{{\tilde{\mu}}}^\xi_{\tilde{x}}(\widetilde{V}^u({f},{\tilde{z}},n,\epsilon/2)) \notag\\
    & \leq\sum_{\tilde{z}\in {\widetilde{F}}}{{\tilde{\mu}}}^\xi_{\tilde{x}}(\widetilde{V}^u({f},\tilde{y}(\tilde{z}),n,\epsilon)).
  \end{align}

  Using \eqref{eq:est2}, \eqref{eq:est3} and Lemma \ref{lem:wellknown} with
  \[
    p_i={{\tilde{\mu}}}^\xi_{\tilde{x}}(\widetilde{V}^u({f},\tilde{y}(\tilde{z}),n,\epsilon))\text{ and }a_i=(\widetilde{S}_n\tilde{\varphi})(\tilde{y}(\tilde{z})),
  \]
  we have
  \begin{align*}
    &\sum_{\tilde{z}\in {\widetilde{F}}}{{\tilde{\mu}}}^\xi_{\tilde{x}}(\widetilde{V}^u({f},\tilde{y}(\tilde{z}),n,\epsilon))\left(n\left(\int_{ M^f}\tilde{\varphi} \mathrm{d}{{\tilde{\mu}}}-\rho\right)+n(h_{{\mu}}({f},\xi)-\rho)\right)\notag\\
\leq &\sum_{\tilde{z}\in {\widetilde{F}}}{{\tilde{\mu}}}^\xi_{\tilde{x}}(\widetilde{V}^u({f},\tilde{y}(\tilde{z}),n,\epsilon))\Big((\widetilde{S}_n\tilde{\varphi})(y(z))-\log {{\tilde{\mu}}}^\xi_{\tilde{x}}(\widetilde{V}^u({f},\tilde{y}(\tilde{z}),n,\epsilon))\Big)\notag\\
\leq &\left(\sum_{\tilde{z}\in {\widetilde{F}}}{{\tilde{\mu}}}^\xi_{\tilde{x}}(\widetilde{V}^u({f},\tilde{y}(\tilde{z}),n,\epsilon))\right)\left(\log \sum_{\tilde{z}\in {\widetilde{F}}}\exp((\widetilde{S}_n\tilde{\varphi})(\tilde{y}(\tilde{z})))-\right. \\
\qquad& \log \left.\sum_{\tilde{z}\in {\widetilde{F}}}{{\tilde{\mu}}}^\xi_{\tilde{x}}(\widetilde{V}^u({f},\tilde{y}(\tilde{z}),n,\epsilon))\right).
  \end{align*}
  Combining  \eqref{eq:est3.5},
  \begin{align}\label{eq:est4}
    & n\Big(\int_{ M^f}\tilde{\varphi} \mathrm{d}{{\tilde{\mu}}}-\rho\Big)+n(h_{{\mu}}({f},\xi)-\rho) \notag\\
    \leq &\log \sum_{\tilde{z}\in {\widetilde{F}}}\exp((\widetilde{S}_n\tilde{\varphi})(\tilde{y}(\tilde{z})))-\log \sum_{\tilde{z}\in {\widetilde{F}}}{{\tilde{\mu}}}^\xi_{\tilde{x}}(\widetilde{V}^u({f},\tilde{y}(\tilde{z}),n,\epsilon))\notag\\
    \leq & \log\sum_{\tilde{z}\in {\widetilde{F}}}\exp((\widetilde{S}_n\tilde{\varphi})(\tilde{y}(\tilde{z})))-\log(1-\rho).
  \end{align}
  Let ${\Delta_{\epsilon}:=\sup\{|\tilde{\varphi}(\tilde{x})-\tilde{\varphi}({\tilde{y}})|\colon d(\Pi(\tilde{x}),\Pi(\tilde{y}))\leq\epsilon\}}$. For any $\tilde{z}\in {\widetilde{F}}$, we have
  \[
    \exp((\widetilde{S}_n\tilde{\varphi})(\tilde{y}(\tilde{z})))\leq\exp((\widetilde{S}_n\tilde{\varphi})({\tilde{z}})+n\Delta_{\epsilon}).
  \]
  Dividing by $n$ and taking the $\limsup$ on both sides of \eqref{eq:est4}, we have
  \[
    \int_{ M^f}\tilde{\varphi} \mathrm{d}{{\tilde{\mu}}}+h_{{\mu}}({f},\xi)-2\rho\leq\limsup_{n\to\infty}\frac{1}{n}\log\sum_{\tilde{z}\in {\widetilde{F}}}\exp((\widetilde{S}_n\tilde{\varphi})({\tilde{z}}))+\Delta_{\epsilon}.
  \]
  We can choose a sequence $\{{\widetilde{F}}_n\}$ of such ${\widetilde{F}}$ such that
  \[
    \limsup_{n\to\infty}\frac{1}{n}\log\sum_{\tilde{z}\in {\widetilde{F}}_n}\exp((\widetilde{S}_n\tilde{\varphi})({\tilde{z}}))\leq {\widetilde{P}^u({\tau},\tilde{\varphi},\delta)}.
  \]
    Since $\rho$ is arbitrary, and $\Delta_{\epsilon}\to 0$ as $\epsilon\to0$, we have
  \[
     \int_{ M^f}\tilde{\varphi} \mathrm{d}{{\tilde{\mu}}}+h_{{\mu}}({f},\xi)\leq \widetilde{P}^u({\tau},\tilde{\varphi},\delta),
  \]
  which implies what we need.
\item In this part, we prove that
  \[
    \sup_{\mu\in\mathcal{M}(f)}\left\{h^u_\mu(f)+\int_{M}\varphi\mathrm{d}{\mu}\right\}=P^u(f,\varphi),
  \]
which completes the proof of Theorem \ref{thm:vp}. In fact, we only need to prove that for any $\rho>0$, there exists ${\tilde{\mu}}\in\mathcal{M} ({\tau})$ such that $h^u_{{\mu}}({f})+\int_{ M^f}\tilde{\varphi} \mathrm{d}{{\tilde{\mu}}}\geq P^u({f},{\varphi})-\rho$.

  Given $\delta>0$, we can choose ${\tilde{x}_0}\in M^f$ such that
  \[
    \widetilde{P}^u({\tau},\tilde{\varphi},{\tilde{x}_0},\delta)\geq \widetilde{P}^u({\tau},\tilde{\varphi},\delta)-\rho.
  \]

Take $\epsilon>0$ small enough. Then let $\widetilde{E}_n$ be an $(n,\epsilon)$ ${W}^u$-separated set of $\overline{\widetilde{W}^u({\tilde{x}_0},\delta)}$ such that
  \[
    \log\sum_{\tilde{y}\in \widetilde{E}_n}\exp((\widetilde{S}_n\tilde{\varphi})(\tilde{y}))\geq\log  \widetilde{P}^u({\tau},\tilde{\varphi},{\tilde{x}_0},\delta,n,\epsilon)-1.
  \]

  Then we construct measures ${\tilde\nu}_{n}$ as follows:
  \[
    {\tilde{\nu}_n}:=\frac{\sum_{\tilde{y}\in \widetilde{E}_n}\exp((\widetilde{S}_n\tilde{\varphi})(\tilde{y}))\tilde{\delta}_{\tilde{y}}}{\sum_{\tilde{z}\in \widetilde{E}_n}\exp((\widetilde{S}_n\tilde{\varphi})(\tilde{z}))},
  \]
  where $\tilde{\delta}_\cdot$ denotes a Dirac measure.
  Let
  \[
    {{\tilde{\mu}}}_{n}=\frac{1}{n}\sum_{i=0}^{n-1}\tau^i{\tilde\nu}_{n}.
  \]
  Then there exists a subsequence $\{n_i\}$ such that
  \[
    \lim_{i\to\infty}{{\tilde{\mu}}}_{n_i}={{\tilde{\mu}}}.
  \]
  It is easy to check that ${\tilde{\mu}}\in\mathcal{M} ({\tau})$.

  We can choose a partition $\eta\in\mathcal{P}^u( M^f)$ such that $\overline{\widetilde{W}^u({\tilde{x}_0},\delta)}\subset \eta_{ }({\tilde{x}_0})$ (by shrinking $\delta$ if necessary). Then choose a finite partition $\alpha$ of $ M^f$ with sufficiently small diameter such that ${{{\mu}}}(\partial \Pi(\alpha))=0$ and suppose that $\alpha$ contains $K$ elements. {Let $\alpha^u$ denote the corresponding measurable partition in $\mathcal{P}^u( M^f)$ constructed via $\alpha$.}

  Fix $q$, $n\in\mathbb{N}$ with $1<q\leq n-1$. Put $a(j)=\left[\frac{n-j}{q}\right]$, $j=0,1,\cdots,q-1$, where we denote by $[a]$ the integer part of $a$. Then
    \[\bigvee_{u=0}^{n-1}\tau^{-i}\alpha=\bigvee_{r=0}^{a(j)-1}\tau^{-(rq+j)}\alpha_0^{q-1}\vee\bigvee_{t\in T_j}\tau^{-t}\alpha,
  \]
  where $T_j=\{0,1,\cdots,j-1\}\cup\{j+aq(j),\cdots,n-1\}$. Note that $\mathrm{Card\ } T_j\leq 2q$. Moreover, we require that $\mathrm{diam}(\alpha)\ll\epsilon$. Then
  \begin{align*}
      & \log\sum_{\tilde{y}\in \widetilde{E}_n}\exp((\widetilde{S}_n\tilde{\varphi})(\tilde{y})) \\
     =& \sum_{\tilde{y}\in \widetilde{E}_n}{\tilde{\nu}_n}(\{\tilde{y}\})\Big(-\log{\tilde{\nu}_n}(\{\tilde{y}\})+(\widetilde{S}_n\tilde{\varphi})(\tilde{y})\Big) \\
     =& H_{\tilde\nu_n}(\alpha^{n-1}_0|\eta)+\int_{ M^f}(\widetilde{S}_n\tilde{\varphi})\mathrm{d}{\tilde\nu}_n.
  \end{align*}
 Then following the same calculation in  \cite{HuWuZhu2017}, we have that
  \begin{align*}
         & \log\sum_{\tilde{y}\in \widetilde{E}_n}\exp((\widetilde{S}_n\tilde{\varphi})(\tilde{y})) \\
    %\leq & \sum_{t\in T_j}H_{{\nu}_n}(\tau^{-t}\alpha|\eta)+H_{\tau^j{\nu}_n}(\alpha^{q-1}_0|\tau^j\eta) \\
    %+    &\sum_{r=1}^{a(j)-1}H_{\tau^{rq+j}{\nu}_n}(\alpha^{q-1}_0|\tau\alpha^u)+\int_{ M^f}(\widetilde{S}_n\tilde{\varphi})\mathrm{d}{\nu}_n\\
     %\leq & \sum_{t\in S_j}\int_\Omega H_{({\nu}}_n)_\omega}(\tau^{-t}\alpha|\eta)d\mathbf{P}(\omega)+H_{\tau^j{\nu}}_n}(\alpha^{q-1}_0|\tau^j\eta) \\
    %+    &\sum_{r=1}^{a(j)-1}H_{\tau^{rq+j}{\nu}}_n}(\alpha^{q-1}_0|\tau\alpha^u)+\int_{ M^f}(\widetilde{S}_n\tilde{\varphi})\mathrm{d}{\nu}}_n\\
    %\leq & \sum_{t\in S_j}\int_{\Omega}\log K(\omega)d\mathbf{P}(\omega)+H_{\tau^j{\nu}}_n}(\alpha^{q-1}_0|\tau^j\eta) \\
    %+    &\sum_{r=1}^{a(j)-1}H_{\tau^{rq+j}{\nu}}_n}(\alpha^{q-1}_0|\tau\alpha^u)+\int_{ M^f}(\widetilde{S}_n\tilde{\varphi})\mathrm{d}{\nu}}_n\\
    \leq & 2q\log K+H_{\tau^j{\tilde\nu}_{n}}(\alpha^{q-1}_0|\tau^j\eta) \\
    +    &\sum_{r=1}^{a(j)-1}H_{\tau^{rq+j}{\tilde\nu}_{n}}(\alpha^{q-1}_0|\tau\alpha^u)+\int_{ M^f}(\widetilde{S}_n\tilde{\varphi})\mathrm{d}{\tilde\nu}_{n}.
  \end{align*}

 Summing the inequality above over $j$ from $0$ to $q-1$ and dividing by $n$, by Lemma \ref{lem:affine} we have
  \begin{align}\label{eq:keyest}
         & \frac{q}{n}\log\sum_{\tilde{y}\in \widetilde{E}_n}\exp((\widetilde{S}_n\tilde{\varphi})(\tilde{y})) \notag \\
    \leq & \frac{2q^2}{n}\log K+\frac{1}{n}\sum_{j=0}^{q-1}H_{\tau^j{\tilde\nu}_{n}}(\alpha^{q-1}_0|\tau^j\eta) \notag  \\
    +    & H_{{{\tilde{\mu}}}_n}(\alpha^{q-1}_0|\tau\alpha^u)+q\int_{ M^f}\tilde{\varphi} \mathrm{d}{{\tilde{\mu}}}_n.
  \end{align}

  Then we can choose a sequence $\{n_k\}$ such that
  \begin{enumerate}[label=(\roman*)]
    \item ${{\tilde{\mu}}}_{n_k}\to{{\tilde{\mu}}}$ as $k\to\infty$;
    \item the following equality holds
    \begin{align*}
     & \lim_{k\to\infty}\frac{1}{n_k}\log \widetilde{P}^u({\tau},\tilde{\varphi},{\tilde{x}_0},\delta,n_k,\epsilon) \\
    =& \limsup_{n\to\infty}\frac{1}{n}\log \widetilde{P}^u({\tau},\tilde{\varphi},{\tilde{x}_0},\delta,n,\epsilon);
    \end{align*}
    \item ${\tilde\nu}_{n_k}\to \tilde\nu$ as $k\to\infty$ for some measure $\tilde\nu$ on $ M^f$.
  \end{enumerate}

  Since ${{\tilde{\mu}}}(\partial \Pi(\alpha))=0$, by Lemma \ref{lem:semiconti},
  \[
    \limsup_{k\to\infty}H_{{{\tilde{\mu}}}_{n_k}}(\alpha_0^{q-1}|\tau\alpha^u)\leq H_{{{\tilde{\mu}}}}(\alpha_0^{q-1}|\tau\alpha^u).
  \]
  As $\tilde\nu_{n}$ is supported on $\overline{\widetilde{W}^u({\tilde{x}_0},\delta)}$, for each $j=0,\cdots,q-1$, we can choose $\alpha, \beta_n \in \mathcal{P}(M^f)$ such that $\beta_1<\beta_2<\cdots<\beta_n<\cdots$ and {$\mathcal{B}(\beta_n)\nearrow\mathcal{B}(\tau^j\eta)$}, and moreover, $(\Pi\tau^j{\tilde\nu})(\partial(\Pi(\alpha_0^{q-1})))=0$, $(\Pi\tau^j{\nu})(\partial(\Pi(\beta_n)_0^{q-1}))=0$. Then applying Lemma \ref{lem:semiconti} we have
  \[
    \limsup_{k\to\infty}\frac{1}{n_k}\sum_{j=0}^{q-1}H_{\tau^j{\tilde\nu}_{n_k}}(\alpha_0^{q-1}|\tau^j\eta)\leq \limsup_{k\to\infty}\frac{1}{n_k}\sum_{j=0}^{q-1}H_{\tau^j{\tilde\nu}}(\alpha_0^{q-1}|\tau^j\eta)=0.
  \]
  Thus replacing $n$ by $n_k$ in \eqref{eq:keyest} and letting $k\to\infty$, by the above claim and discussions, we get
  \begin{align*}
        & q\limsup_{n\to\infty}\frac{1}{n}\log \widetilde{P}^u({\tau},\tilde{\varphi},{\tilde{x}_0},\delta,n,\epsilon) \\
    \leq& H_{{{\tilde{\mu}}}}(\alpha^{q-1}_0|\tau\alpha^u)+q\int_{ M^f}\tilde{\varphi}\mathrm{d}{{\tilde{\mu}}}.
  \end{align*}
  By Theorem \ref{thm:localvsfinite},
  \begin{align*}
        & \limsup_{n\to\infty}\frac{1}{n}\log \widetilde{P}^u({\tau},\tilde{\varphi},{\tilde{x}_0},\delta,n,\epsilon) \\
    \leq& \lim_{q\to\infty}\frac{1}{q}H_{{{\tilde{\mu}}}}(\alpha^{q-1}_0|\tau\alpha^u)+\int_{ M^f}\tilde{\varphi}\mathrm{d}{{\tilde{\mu}}} \\
       =   & h^u_{{\mu}}({f})+\int_{ M^f}\tilde{\varphi}\mathrm{d}{{\tilde{\mu}}}.
  \end{align*}
  Let $\epsilon\to 0$, we have $\widetilde{P}^u({\tau},\tilde{\varphi},{\tilde{x}_0},\delta)\leq  h^u_{{\mu}}({f})+\int_{ M^f}\tilde{\varphi}\mathrm{d}{{\tilde{\mu}}}$. Recall that $P^u({f},{\varphi})=\widetilde{P}^u({\tau},\tilde{\varphi}, \delta)\leq \widetilde{P}^u({\tau},\tilde{\varphi},{\tilde{x}_0},\delta)+\rho$. The proof of Theorem \ref{thm:vp} is complete.
\end{enumerate}
\end{potC}
\subsection{\texorpdfstring{$u$-equilibrium states for endomorphisms}{u-equilibrium states for endomorphisms}}
In this subsection, we introduce the notion of $u$-equilibrium state and list some results concerning it, whose proofs are similar to those in  \cite{HuWuZhu2017}.
Let $\varphi\in C(M)$.
\begin{mdef}\rm\label{def:uequi}
  $\mu\in\mathcal{M}({f})$ is said to be a $u$-equilibrium state for $\varphi$, if it satisfies
  \[
  h^u_\mu({f})+\int_{ M^f}\varphi\mathrm{d}{\mu}=P^u({f},\varphi).
  \]
\end{mdef}
We denote by $\mathcal{M}^u({f},\varphi)$ the set of all $u$-equilibrium states for $\varphi$.
\begin{prop}\label{prop:equilibrium}
  Let $\varphi\in C(M)$, then we have the following properties related with $u$-equilibrium states.
  \begin{enumerate}[label=(\roman*)]
    \item\label{equi1} $\mathcal{M}^u({f},\varphi)$ is non-empty, and it is convex, in particular, the measure of maximal unstable metric entropy always exists;
    \item the extreme points of $\mathcal{M}^u({f},\varphi)$ are precisely ergodic members of $\mathcal{M}^u({f},\varphi)$;
    \item $\mathcal{M}^u({f},\varphi)$ is compact and has an {ergodic $u$-equilibrium state};
    \item assume $\varphi$, $\psi\in C(M)$ are cohomologous, i.e. $\varphi=\psi+\sigma-\sigma\circ\tau-c$ for some $c\in\mathbb{R}$ and $\sigma\in C(M)$. Then $\varphi$ and $\psi$ have the same $u$-equilibrium states, and
        \[
          P^u({f},\varphi)=P^u({f},\psi)-c.
        \]
  \end{enumerate}
\end{prop}

\section*{Acknowledgements}

X. Wang and Y. Zhu are supported by  NSFC (No: 11771118, 11801336), W. Wu is supported by  NSFC (No: 11701559). The first author is also supported by  the Innovation Fund Designated for Graduate Students of Hebei Province (No: CXZZBS2018101) and China Scholarship Council (CSC).

The authors would like to thank the referees for the detailed review and very valuable suggestions, which led to improvements of the paper.

%\section*{References}

{xswang@hebtu.edu.cn (Xinsheng Wang)}

{wuweisheng@cau.edu.cn (Weisheng Wu)}

{yjzhu@xmu.edu.cn (Yujun Zhu)}
\end{document}